\documentclass[final,hidelinks,onefignum,onetabnum]{siamart220329}

\usepackage{lipsum}
\usepackage{amsfonts}
\usepackage{amsmath}
\usepackage{amssymb}
\usepackage{graphicx}
\usepackage{epstopdf}
\usepackage{algorithmic}
\usepackage{amsopn}
\usepackage{tikz}      
\usepackage{pgfplots}   
\pgfplotsset{compat=1.18}
\usepackage{nicefrac}
\usepackage{comment}
\usepackage{subcaption}
\usepackage{hyperref}
\usepackage{fontawesome}
\usepackage{bm}
\usepackage{enumitem}
\usepackage{float}  
\usepackage{mathdots} 
\usepackage{array}
\usepackage{booktabs}
\usepackage{yfonts}
\usepackage{arydshln}
\usepackage{bm}
\usetikzlibrary{arrows, arrows.meta, positioning, fit, calc}

\definecolor{istblue}{RGB}{0,65,145}

\newsiamremark{rem}{Remark}
\newsiamthm{claim}{Claim}
\newsiamremark{hypothesis}{Hypothesis}
\crefname{hypothesis}{Hypothesis}{Hypotheses}
\newsiamremark{exmp}{Example}
\crefname{exmp}{Example}{Example}
\newsiamthm{problem}{Problem}
\crefname{problem}{Problem}{Problem}
\newsiamthm{assum}{Assumption}
\crefname{assum}{Assumption}{Assumption}
\newsiamthm{prop}{Proposition}
\crefname{prop}{Proposition}{Proposition}
\newsiamthm{cor}{Corollary}
\crefname{cor}{Corollary}{Corollary}

\newcommand{\norm}[1]{\| #1 \|}
\newcommand{\R}{\mathbb{R}}

\newcommand{\Smlt}{\mathcal{S}_\theta(m,L)}
\newcommand{\Vtnu}{\mathcal{V}_{\Theta,\nu}}

\allowdisplaybreaks

\definecolor{mycolor3}{rgb}{0.0329, 0.4430, 0.8720}
\definecolor{mycolor4}{rgb}{0.1959, 0.2645, 0.7279}
\definecolor{mycolor2}{rgb}{0.3482, 0.7424, 0.5473}
\definecolor{mycolor1}{rgb}{0.6834, 0.7435, 0.4044}

\headers{Analysis of Varying Algorithms with Robust Control.}{Fabian Jakob and Andrea Iannelli}

\title{A Linear Parameter-Varying framework for the analysis of time-varying optimization algorithms 
}

\author{Fabian Jakob\thanks{Institute for Systems Theory and Automatic Control, University of Stuttgart, Germany (\email{fabian.jakob@ist.uni-stuttgart.de}, \email{andrea.iannelli@ist.uni-stuttgart.de}).}
	\and Andrea Iannelli\footnotemark[1]}
	
\ifpdf
\hypersetup{
	pdftitle={A Linear Parameter-Varying framework for the analysis of time-varying optimization algorithms },
	pdfauthor={F. Jakob, A. Iannelli}
}
\fi

\ifpdf
\DeclareGraphicsExtensions{.eps,.pdf,.png,.jpg}
\else
\DeclareGraphicsExtensions{.eps}
\fi

\makeatletter
\def\cref@currentlabel{[][]{}}
\makeatother

\begin{document}

\maketitle

\begin{abstract}
    In this paper we propose a framework to analyze iterative first-order optimization algorithms for time-varying convex optimization. We assume that the temporal variability is caused by a time-varying parameter entering the objective, which can be measured at the time of decision but whose future values are unknown. We consider the case of strongly convex objective functions with Lipschitz continuous gradients. We model the algorithms as discrete-time linear parameter varying (LPV) systems in feedback with the time-varying gradient. We leverage the approach of analyzing algorithms as uncertain control interconnections with integral quadratic constraints (IQCs) and generalize that framework to the time-varying case. We propose novel IQCs that are capable of capturing the behavior of time-varying nonlinearities and leverage techniques from the LPV literature to establish novel bounds on the tracking error. Quantitative bounds can be computed by solving a semi-definite program and can be interpreted as an input-to-state stability result with respect to a disturbance signal which increases with the temporal variability of the problem. As a departure from results in this research area, our bounds introduce a dependence on different additional measures of temporal variations, such as the rate of change of function or gradient value.
    We exemplify our main results with numerical experiments that showcase how our analysis framework is able to capture convergence rates of different first-order algorithms for time-varying optimization through the choice of IQC and rate bounds.
\end{abstract}



\section{Introduction}\label{sec:intro}

Time-varying convex optimization is an emerging subfield of convex optimization in which the optimization problem exhibits some temporal variability in the objective function, the constraints, or both \cite{dallAnese_informationStreams}. The development of algorithms that are able to accurately track the time-varying solution of such problems can be quite impactful, as application domains include, for instance, power grid systems \cite{powerFlowPursuit}, mobile robotics \cite{pappas_mobileRobotics}, congestion control \cite{low_congestion}, or traffic systems \cite{bianchin_internalmodel}.
Recent research efforts resulted in many algorithms that are able to asymptotically track optimal solutions if a model of temporal variability is available \cite{fazylab, davydov, bastianello_internalModel, bianchin_internalmodel}. However, in many practical scenarios, having this knowledge might be unrealistic or even impossible. Often, the remedy is then to resort to some simple first-order algorithm such as gradient descent or primal-dual methods, cf. \cite{Bernstein,colombino2018,bianchinFeedbackOpt}. These algorithms can be studied quite conveniently, as bounds on asymptotic tracking can readily be derived; see \cite{popkov} for one of the first works. On the other hand, the effect of momentum, by now well understood in the static case \cite{nesterov,triple_momentum}, has not received much attention for time-varying problems. Performance degradation, e.g. shown empirically in \cite{dallAnese_informationStreams} for Nesterov's fast gradient method, and possibly tedious convergence proofs are important motivations for developing a framework to understand and analyze general first-order algorithms for time-varying optimization.

In recent years, a powerful framework based on robust control methods has been developed to study general first-order methods for static (or time-invariant) convex optimization, which is able to \emph{systematically} analyze their convergence and robustness \cite{lessard, lessardSeilerSynthesis, jovanovic_robustness, michalowski, scherer_ebenbauer, lessardDissipativity}. The basic idea is to model an algorithm as a discrete-time linear time-invariant (LTI) system in feedback with the gradient, and leveraging \emph{integral quadratic constraints} (IQCs) \cite{megretzki} to model the input-output behavior of the (unknown) nonlinear part. This ultimately yields an absolute stability problem which has been thoroughly studied in the robust control literature \cite{zames_falb}. Starting from this powerful modelling viewpoint, further extensions have been proposed, e.g. the incorporation of a performance channel to analyze the effect of gradient noise \cite{jovanovic_robustness} and the development of convex conditions for synthesis of algorithms \cite{lessardSeilerSynthesis, scherer_ebenbauer, gramlich}.

Motivated by its advantages and strengths, this work generalizes this approach to time-varying convex optimization problems and develops a framework that systematically provides tracking guarantees and quantifiable convergence bounds for general first-order algorithms. We assume the temporal variability of the optimization problem is caused by a time-varying parameter entering the objective function and propose to model the algorithms as interconnections of \emph{linear parameter-varying} (LPV) systems in feedback with the time-varying gradient. Consequently, the convergence proofs and IQCs to describe varying nonlinearities have to be extended or newly developed. To address the challenges arising from this new problem setting, we develop a holistic analysis blending tools from time-varying optimization, IQCs, and LPV theory.

\textbf{Related works.}\,\, A general overview of algorithms for time-varying convex optimization can be found in \cite{simonetto_timeStructuredReview}, 
both for methods that use information of temporal variability, i.e. through the use of a known time-derivative of parameter or gradient, and those that do not. In \cite{bianchin_internalmodel} it is shown that leveraging knowledge of temporal variability is both necessary and sufficient for exact asymptotical tracking of the time-varying minimizer. Convergence results of algorithms that do not have such knowledge can only be established to neighborhoods of the optimal trajectory and typically arise in the form of input-to-state stability (ISS) conditions \cite{sontag, jiangWang_iss}, where some signal quantifying the temporal variability acts as the disturbance \cite{Bernstein,colombino2018, bianchinFeedbackOpt, tabuada}.
In \cite{dallAnese_krasovsii, simonetto_averagedOperators} tracking guarantees are provided for more general first-order running algorithms in the form of time-varying averaged operators. Those algorithms have the advantage of being able to tackle non-smooth optimization problems on the one hand, however, the results rely on a boundedness assumption on the operator; a condition that is typically only met in constrained optimization problems. Moreover, the proposed bounds can be conservative in application, as the tracking bound depends on the diameter of the feasible set.

The extension of the IQC framework for first-order algorithms to time-varying algorithms is mentioned as an outlook in the seminal work \cite{lessard}. A few attempts to pursue this direction have been made in the literature. Gradient descent with varying step size has been framed in \cite{padmanabhan} as an LPV system, but the analysis essentially treats the stepsize only as a parametric uncertainty rather than leveraging LPV techniques. In \cite{fazylab} algorithms with time-varying parameters are considered and handled through time-varying Lyapunov functions, but subsequent analyses rely on case-by-case inspections. Moreover, both works consider time-varying algorithms for static optimization only. In \cite{simonettoMassionini} LPV tools are used to synthesize an optimal Kalman filter for time-varying problems, however, the scheduling parameter explicitly depends on the temporal variability of the problem, which is an assumption we do not make in this work. Moreover, unlike in this work, the problem is not cast as an LPV-IQC augmented plant, so that an unified analysis of the proposed interconnection is not possible. Despite time-varying costs being so far disregarded in the IQC framework, there do exist works on analyzing uncertain LPV systems with IQCs \cite{pfifer, takarics, wang_lpv_iqc_syn}, albeit in continuous time. The application to time-varying first-order algorithms requires an extension to discrete time under consideration of changing set points.
 
\textbf{Contributions.}\,\, In this work, we propose a novel and general approach to model first-order algorithms for smooth time-varying optimization problems by casting them as LPV systems in feedback with a time-varying gradient. Our algorithm formulation encompasses both algorithms that perform one or multiple iterations until the objective changes. Our setting considers parameter-varying IQCs and LPV systems whose set points fulfill the time-varying optimality condition, thus tailoring the setup to the use case of time-varying convex optimization. To this aim, we propose a 
more general type of IQCs that can characterize the input-output behavior of time-varying gradients of strongly convex and Lipschitz-smooth cost functions, up to some interpretable residual term that vanishes for static problems. Crucially, the Lipschitz and convexity constants, as well as all algorithm's parameters, are all allowed to be time-varying. We show how the analysis results allow obtaining computable convergence bounds for this more general setting and provide a new systematic proof strategy for general first-order time-varying optimization algorithms.

\textbf{Outline.}\,\, The paper is structured as follows. The problem statement is given in \cref{sec:problem_setting} and the algorithm formulation is introduced in \cref{sec:algorithm_formulation}.
In \cref{sec:iqcs} we derive a new type of IQCs to handle time-varying slope restricted nonlinearities. Based on them, we establish our tracking bounds in \cref{sec:tracking}.Within this section we will highlight the difference between the bounds that can be established with pointwise quadratic constraints and dynamic IQCs. Our bounds can be computed efficiently by solving a semidefinite program (SDP) in the form of a linear matrix inequality (LMI). We additionally present case studies in \cref{sec:case_studies} investigating the influence of algorithmic structures and rate bounds on the convergence rate. 
Finally, we conclude the paper in \cref{sec:conclusion}.

\textbf{Notation.}\,\,
Let $\mathbb{N}_0$ denote the set of natural numbers including $0$ and let $\mathbb{S}^n \subset \R^{n\times n}$ denote the set of real-valued symmetric matrices. 
We denote by $I_n$ the $n\times n$ identity matrix, by $1_n$ a column vector of ones of length $n$, and by $\mathbb{I}_n$ the index set of integers from $1$ to $n$. If $P\in \mathbb{S}^n $ is positive (semi-)definite, we write $P\succ 0$ ($\succeq$) and denote the weighted square norm as $\| x \|_P^2 = x^\top P x$. We omit the index for the Euclidean norm, i.e., for $P=I_n$. For a matrix $U$, $\| U \|_2$ denotes the induced 2-norm. For two column vectors $x$ and $y$ we will write $\mathrm{vec}(x,y):=\begin{bmatrix}x^\top & y^\top\end{bmatrix}^\top$. The Kronecker product is denoted by $\otimes$. 

Given some initialization $\xi_0$, a discrete-time linear dynamical system with internal state $\xi_k$, input sequence $\{u_k\}_{k\geq 0}$ and output sequence $\{y_k\}_{k\geq 0}$ is given by 
\begin{equation*}
    \xi_{k+1} = A \xi_k + B u_k, \quad y_k = C \xi_k + D u_k, \quad k\geq 0
\end{equation*}
where the matrices have compatible dimensions and may be parameter-varying, depending on the context. We compactly express this $u \mapsto y$ mapping as
\begin{equation*}
y_k=\left[
    \begin{array}{c|c}
    A & B \\ \hline
    C & D
    \end{array}
\right]_{\xi_0} u_k.
\end{equation*}
We omit the index $\xi_0$ if $A,B,C=0$ (state-less feedthrough mapping $y_k = Du_k$).

\section{Problem Setting}\label{sec:problem_setting}

We consider the unconstrained time-varying optimization problem
\begin{equation}\label{eq:min_fx_varying}
	x^\star_k := x^\star(\theta_k) = \arg \min_{x \in \R^d} f(x, \theta_k)
\end{equation}
with time index $k\in \mathbb{N}_0$ and objective function $f: \R^d \times \Theta \rightarrow \R$. 
The time-variation is driven by a parameter sequence $\theta: \mathbb{N}_0 \rightarrow \Theta$, which belongs to some compact domain~$\Theta \subset \R^{n_\theta}$. 
We seek to take a decision $x_{k+1}$ using some first-order algorithm and assume that the parameter $\theta_k$ is measurable at each timestep but no knowledge on future $\theta_{k+t}, t\geq 1$ is available, yielding an \emph{information-limited} setup \cite{dallAnese_informationStreams,Bernstein,simonetto_averagedOperators}.
We make the following assumption on the objective function $f$.
\begin{assum}
The objective $f(x,\theta)$ is $m(\theta)$-strongly convex and $L(\theta)$-smooth in $x$, where $m, L : \Theta \to \mathbb{R}$ are given continuous functions satisfying $0 < \underline{m} < m(\theta) < L(\theta) < \overline{L} < \infty$ for all $\theta \in \Theta$ and some scalars $\underline{m}$, $\overline{L}$. In addition, $\nabla_x f(x,\theta)$ is continuous in $\theta$ for every $x\in\mathbb{R}^d$.
\end{assum}
The set of such objective functions will be compactly expressed as $\Smlt$.
Such regularity assumptions are standard in time-varying convex optimization (cf. \cite{bianchinFeedbackOpt,popkov,simonetto_timeStructuredReview}), which on the one hand ensure existence and uniqueness of the minimizer $x^\star_k$ at each $k$, and on the other enable the establishment of linear convergence-type results.
We define the tracking error as $\|x_k -x_k^\star\|$, where $\{x_k\}_{k\geq 0}$ is the sequence of decisions of some algorithm and $\{x_k^\star\}_{k\geq 0}$ the optimal trajectory solving \cref{eq:min_fx_varying}.

\subsubsection*{The class of first-order algorithms}

We consider the class of first-order algorithms that can be written as discrete-time linear dynamical systems
\begin{subequations}\label{eq:algorithm_lpv}
    \begin{align}
    \label{eq:system_lpv}
        \begin{aligned}
            \xi_{k+1} & = A(\theta_k) \xi_k + B(\theta_k) u_k \\
            y_k       & = C(\theta_k) \xi_k + D(\theta_k) u_k,
        \end{aligned} 
    \end{align}
with internal state $\xi_k \in \mathbb{R}^{n_\xi}$, output $y_k \in \R^{n_y}$, and an input $u_k \in \R^{n_y}$ that encodes first-order information of \cref{eq:min_fx_varying} through some nonlinear feedback 
    \begin{equation}\label{eq:oracle}
        u_k = \varphi_k(y_k).
    \end{equation}

The algorithm's decisions are extracted through a state readout equation
\begin{equation}\label{eq:readout}
    x_{k+1} = C_x\, \xi_{k+1}.
\end{equation}

\end{subequations}

In particular, $\varphi$ will include (possibly multiple) evaluations of the gradient $\nabla_x f( \cdot, \theta_k)$, or other first-order oracles in more general settings (e.g., projections for constrained problems). Given some initialization $\xi_0$, \cref{eq:system_lpv,eq:oracle} defines a sequence of states $\{ \xi_k \}_{k\geq 0}$ from which the decisions are extracted. We assume throughout that $A,B,C,D$ are continuous functions of the parameter $\theta$. Conceptually, we will use this formulation to express algorithms where the decision $x_{k+1}$ depends on the current and past iterates $x_k, x_{k-1}, \dots$; available parameter values $\theta_k, \theta_{k-1}, \dots$; as well as on first-order information available at each step. We illustrate this in the next section. Such dynamical system formulations of first-order algorithms have been already proposed in static optimization settings, e.g. in \cite{michalowski,lessardDissipativity,upadhyaya}. The novelty of this work is to propose a formulation that captures also a broad class of first-order methods relevant for time-varying optimization problems. Notably, \cref{eq:system_lpv} yields what in control theory is known as an LPV system \cite{shamma}. Its combination with \cref{eq:oracle} yields a feedback interconnection, see Fig.~\ref{fig:algorithm} for a visualization. 
\begin{figure}[ht!]
    \centering
    % lower arrows + labels are gray now
\begin{tikzpicture}[
  block/.style = {draw, minimum width=0.5cm, minimum height=0.5cm, align=center, fill=white},
  vecarrow/.style = {->, double, line width=0.3pt, >=Latex, double distance=0.8pt},
  vecarrowop/.style = {<-, double, line width=0.3pt, >=Latex, double distance=0.8pt},
  stacklayer/.style = {
    draw, minimum width=0.5cm, minimum height=0.5cm, align=center, fill=white
  },
  ->,
  >=Latex
]

% Gradient block stack
\begin{scope}
  \foreach \i/\j in {2/2, 1/1} {
    \node[stacklayer] (Grad\i) at (-\i*0.1, \j*0.1) {\phantom{$\nabla f_k$}};
  }
  \node[block] (Grad) at (0,0) {$\nabla f_k$};
\end{scope}

% Algorithm block
\node[block, below=2em of Grad] (Algo) {
  $\begin{aligned}
    \xi_{k+1} &= A(\theta_k) \xi_k + B(\theta_k) u_k \\
    y_k &= C(\theta_k) \xi_k + D(\theta_k) u_k
  \end{aligned}$
};

% Arrows (only gradient channel remains)
\draw[vecarrowop] (Grad.east) -- ++(2.5,0) |- node[anchor=west, pos=0.5, yshift=2em] {$y_k$}  ($(Algo.east) + (0,0.0)$);
\draw[vecarrow] (Grad1.west) --  ++(-2.5,0) |- node[anchor=east, pos=0.5, yshift=2em] {$u_k$} ($(Algo.west) + (0,0.0)$);

% % lower arrows made gray
% \draw[->, gray] ($(Algo.west) + (-0.5,-0.3)$) -- node[anchor=east, left, pos=0, text=gray] {Init. $\xi_0$} ($(Algo.west) + (0,-0.3)$);

% \draw[<-, gray] ($(Algo.east) + (0.5,-0.3)$) -- node[anchor=east, right, pos=0, text=gray] {$x_k = C_x \xi_k$} ($(Algo.east) + (0,-0.3)$);

% Dashed box around gradient block only
\node[draw=black, dashed, inner xsep=0.5cm, inner ysep=6pt, fit=(Grad2)(Grad1)(Grad)] (PhiBox) {};
\node[anchor=west, yshift=0.4cm] at (PhiBox.east) {$\varphi_k$};

\end{tikzpicture}
    \caption{The first-order algorithm \cref{eq:algorithm_lpv} as a feedback interconnection of an LPV system and an oracle $\varphi$, which contains (multiple) evaluations of the gradient.}
    \label{fig:algorithm}
\end{figure}
We will characterize the parametrization such that its fixed-point encodes a solution to \cref{eq:min_fx_varying}, and derive convergence guarantees via stability properties of \cref{eq:algorithm_lpv}.

\subsubsection*{Measures of temporal variability}

To capture the temporal evolution of \cref{eq:min_fx_varying}, we will introduce different variability measures. In particular, we introduce the minimizer, function, and gradient change as
\begin{subequations}\label{eq:variational}
\begin{samepage}
\begin{align}
    \label{eq:variational:x}
    \Delta x^\star_k &:= x_k^\star - x_{k+1}^\star, \\
    \label{eq:variational:f}
    \Delta f_k (\cdot) &:= f(\cdot, \theta_k) - f(\cdot, \theta_{k+1}), \\
    \label{eq:variational:grad}
    \Delta g_k (\cdot) &:= \nabla_x f(\cdot, \theta_k) - \nabla_x f(\cdot, \theta_{k+1}),
\end{align}
\end{samepage}
\end{subequations}
respectively. To model the temporal evolution of the parameter sequence itself, we assume a bounded rate of variation for $\theta_k$, i.e.
\begin{equation}\label{eq:param_rate_bound}
    \underline{\nu} \leq \theta_{k+1} - \theta_{k} \leq \overline{\nu},
\end{equation}
where $\underline{\nu}, \overline{\nu} \in \mathbb{R}^{n_\theta}$ are some known bounds and the inequality has to be interpreted componentwise.

Note that compactness of $\Theta$, together with \cref{eq:param_rate_bound}, implies that at each time $k$, $\theta_k$ and $\Delta \theta_k := \theta_{k+1} - \theta_{k}$ form a set of admissible parameter-variation pairs. We capture this relationship with a set $\mathcal{V}_{\Theta,\nu}$, that is defined as 
\begin{equation}\label{eq:consistent_polytope}
\Vtnu = \{ (\theta, \Delta \theta):  \,\,\, \theta \in \Theta, \quad \theta+\Delta \theta \in \Theta, \quad \underline{\nu} \leq \Delta \theta \leq \overline{\nu}\},
\end{equation}
see Fig.~\ref{fig:consistent_polytope} for an illustration. Such a set excludes large variations $\Delta \theta_k$ near the boundary of $\Theta$, which is a concept that originates in LPV stability analysis \cite{deSouza_consistentPolytope, karimi, amato}.

\begin{figure}[H]
    \centering
    \definecolor{anth80}{RGB}{216, 218, 219}

\begin{tikzpicture}[scale=1.3]
		
    % Pentagon shape
    \filldraw[fill=anth80, draw=black] (-1,0) -- (-1,1) -- (0.5,1) -- (1.5,0) -- (1.5,-0.5) -- (-0.5,-0.5) -- cycle;
    
    % Axes
    \draw[->] (-2,0) -- (2.5,0) node[right] {$\theta$};
    \draw[->] (0,-1) -- (0,1.75) node[pos=0.95, right] {$\Delta \theta$};
    
    % dashed corners
    \draw[dashed] (1.5,0) -- (1.5,1) -- (0.5,1);
    \draw[dashed] (-1,0) -- (-1,-0.5) -- (-0.5,-0.5);
    
    % Labels
    \node at (-1,0) [left] {};
    \node[anchor=south west] at (0,1)  {$\overline{\nu}$};
    \node[anchor=south west] at (1.5,0) {$\theta_{\max}$};
    \node[anchor=north west] at (0,-0.5)  {$\underline{\nu}$};
    \node[anchor=south east] at (-1,0) {$\theta_{\min}$};

\end{tikzpicture}
    \caption{Illustration of the set $\Vtnu$ for a one-dimensional box-constrained parameter.}
    \label{fig:consistent_polytope}
\end{figure}

\subsubsection*{Performance quantification}

To quantify the tracking performance of an algorithm, we consider bounds on the tracking error $\| x_k - x_k^\star \|$.
Related work in time-varying optimization (cf. \cite{Bernstein,colombino2018,bianchinFeedbackOpt,popkov}) have shown that for different first-order algorithms one can obtain exponential tracking bounds of the form
\begin{equation}\label{eq:tracking_bound}
	\| x_k - x_k^\star \| \leq c_1 \rho^k \| x_0 - x_0^\star \| + c_2 \sum_{t=1}^{k} \rho^{k-t} \delta_t,
\end{equation}
with some constants $c_1, c_2 > 0$, decay rate $\rho\in(0,1)$, and $\delta_k = \| \Delta x^\star_{k-1} \|$ as a measure of temporal variability\footnote{Note that some of the mentioned works derive their results in continuous time and for constrained problems by using primal-dual methods on the Lagrangian or projected gradient descent, where the discretization and/or the reduction to the unconstrained case would yield \cref{eq:tracking_bound}.}. From a system theoretic point of view, \cref{eq:tracking_bound} can be interpreted as an input-to-state stability result \cite{jiangWang_iss}, and with a norm bound on $\delta_k$ one can establish practical stability \cite{practical_stability}. In other words, the iterates exponentially converge to a ball around the trajectory of minimizers; and the size of that ball depends on the temporal variability of the problem, where past variations are exponentially forgotten in time.
Such bounds are often the best achievable in information-limited settings, but the choice of algorithm can affect both the rate $\rho$ and the size of the ball.
We aim to provide a systematic algorithm analysis procedure that allows one to retrieve this type of results for generic first-order algorithms by posing the following problem:
\begin{problem}\label{prob:tracking_bound}
    Given problem \cref{eq:min_fx_varying}, parameter rate bounds \cref{eq:param_rate_bound}, and a general first-order algorithm of the form \cref{eq:algorithm_lpv}, derive conditions on the algorithm parametrization $A(\theta), B(\theta)$, $C(\theta)$, $D(\theta)$, such that the tracking error $\|{x}_k - x_k^\star\|$ can be bounded by \cref{eq:tracking_bound}, where $\delta_k$ is some term reflecting the temporal variability of \cref{eq:min_fx_varying}.
\end{problem}
As a result, the solution to \cref{prob:tracking_bound} will yield a computer-aided analysis tool that provides tracking certificates for any algorithm that can be formulated as \cref{eq:algorithm_lpv}.
We will build on the framework of IQC-based algorithm analysis to address \cref{prob:tracking_bound}, and develop novel analysis techniques to account for the time-varying nature of the problem.

\section{Algorithm Formulation}\label{sec:algorithm_formulation}
Our starting point is the interpretation of first-order algorithms as dynamical systems. For notational convenience, we use the convention $\nabla f_k(\cdot) \triangleq \nabla_x f(\cdot, \theta_k)$. To show the generality of the class of algorithms considered here, consider the following example.

\begin{exmp}[Accelerated gradient method]\label{exmp:accelerated_algo}
	Consider the general form of an accelerated gradient method
	\begin{equation}\label{eq:accelerated_gm:pure}
		x_{k+1} = x_k + \beta_k (x_k - x_{k-1}) - \alpha_k \nabla f_k\bigl(x_k + \gamma_k (x_k - x_{k-1})\bigr).
	\end{equation}
	Here, $\alpha_k > 0$ and $\beta_k, \gamma_k \geq 0$ are the (potentially time-varying) stepsize and momentum parameters, respectively. Depending on their choice, we get gradient descent (GD), Nesterov's method (NM), the heavy-ball method (HB), or triple momentum (TM), cf. \cite[Table 1]{michalowski}. Typically, $\alpha_k, \beta_k$ and $\gamma_k$ are computed based on the convexity and smoothness moduli $m(\theta_k)$ and $L(\theta_k)$, and thus, they can be considered functions of the parameter $\theta_k$. Now note that \cref{eq:accelerated_gm:pure} can be equivalently written as
	\begin{align}\label{eq:accelerated_gm:ss}
		\begin{aligned}
		\begin{bmatrix} x_{k+1} \\ x_k \end{bmatrix}
		&=
		\begin{bmatrix} (1 + \beta_k) I_d & -\beta_k I_d \\ I_d & 0 \end{bmatrix}
		\begin{bmatrix} x_k \\ x_{k-1} \end{bmatrix}
		-
		\begin{bmatrix} \alpha_k I_d \\ 0 \end{bmatrix}
		u_k, \quad u_k = \nabla f_k(y_k) \\
		y_k &= \begin{bmatrix} (1 + \gamma_k)I_d & -\gamma_k I_d \end{bmatrix} \begin{bmatrix} x_k \\ x_{k-1} \end{bmatrix}.
		\end{aligned}
	\end{align}
	Observe that \cref{eq:accelerated_gm:ss} is precisely of the form \cref{eq:algorithm_lpv} with state $\xi_k \triangleq \mathrm{vec}(x_k, x_{k-1})$ and oracle $\varphi_k \triangleq \nabla f_k$, where the matrices $A(\theta_k),B(\theta_k),C(\theta_k)$ can be read from \cref{eq:accelerated_gm:ss} and $D=0$. The algorithm decision is obtained with $C_x = \begin{bmatrix} I_d & 0 \end{bmatrix}$.
\end{exmp}

\cref{exmp:accelerated_algo} is frequently used in the classic IQC-based algorithm analysis literature (e.g. \cite{lessard, michalowski, scherer_ebenbauer}). There is one notable difference in the time-varying setting: the changing problem data is reflected in changing algorithm parameters (e.g. a popular stepsize tuning law would result in $\alpha_k = 1/L(\theta_k)$), which motivates the LPV formulation.

Another feature that  is unique in the time-varying settings is illustrated in the next example.

\begin{exmp}[Two-step gradient descent]\label{exmp:two_step}
	Consider a scenario in which at time~$k$ we have enough computational budget to perform two steps of gradient descent before new problem data arrives. The algorithm is
	\begin{align}
		\label{eq:example:m_step}
		\begin{aligned}
		\bar{x}_{k} &= x_k - \alpha_k \nabla f_k(x_k)\\
        x_{k+1} &= \bar{x}_{k} - \alpha_k \nabla f_k(\bar{x}_{k}),
		\end{aligned}
	\end{align}
	or, in condensed form, $x_{k+1} = x_k - \alpha_k \nabla f_k(x_k) - \alpha_k \nabla f_k($\footnotesize$x_k - \alpha_k \nabla f_k(x_k)$\normalsize$)$. Now define $y_k := \mathrm{vec}(x_k, \bar x_k)$. Then, \cref{eq:example:m_step} can be equivalently written as
	\begin{align}
		\begin{aligned}\label{eq:example:m_step:ss}
		x_{k+1} &= x_k + \begin{bmatrix} -\alpha_k I_d & -\alpha_k I_d \end{bmatrix} u_k, \quad u_k = \begin{bmatrix} \nabla f_k(y_k^{1}) \\ \nabla f_k(y_k^{2}) \end{bmatrix}\\
		y_k &= \begin{bmatrix} I_d \\ I_d \end{bmatrix} x_k +  \begin{bmatrix} 0 & 0 \\ -\alpha_k I_d & 0 \end{bmatrix} u_k.
		\end{aligned}
	\end{align}
	This formulation is again precisely of the form \cref{eq:algorithm_lpv} with state $\xi_k \triangleq x_k$, the oracle $\varphi_k(y_k) \triangleq \mathrm{vec}(\nabla f_k(y_k^{1}), \nabla f_k(y_k^{2}))$ and $C_x = I_d$.
\end{exmp}

In the same spirit also generic $p$-step versions ($p\geq 2$) of gradient descent can be formulated. 
\cref{exmp:two_step} highlights an important aspect of time-varying optimization: there are two timescales present, one associated with the evolution of the optimization problem and one in which the algorithm iterations change, see Fig.~\ref{fig:multistep:a} for an illustration. 

\begin{figure}[ht]
	\centering
	\begin{subfigure}[t]{0.47\textwidth}
		\centering
		\begin{tikzpicture}[>=latex, scale=1.2]

% Parameters
\def\nslow{2}    % number of points on slow timescale
\def\nfast{4}    % number of points on fast timescale
\def\L{3}        % total length

% Vertical offset for the fast line
\def\offset{-0.9}

% Draw slow timescale
\draw[thick] (0,0) -- (\L,0);
\foreach \i in {0,...,\nslow} {
  \pgfmathsetmacro\x{(\i/\nslow)*\L}
  \draw[thick] (\x,0.1) -- (\x,-0.1);
  \filldraw[black] (\x,0) circle (1.5pt);
  \node[above] at (\x,0.1) {\small $k=\i$};
}

% Draw fast timescale
\draw[thick] (0,\offset) -- (\L,\offset);
\foreach \i in {0,...,\nfast} {
  \pgfmathsetmacro\x{(\i/\nfast)*\L}
  \draw[thick] (\x,\offset+0.1) -- (\x,\offset-0.1);
  \filldraw[istblue] (\x,\offset) circle (1.2pt);
  % Labels
  \ifnum\i=0
    \node[below] at (\x,\offset-0.1) {\footnotesize $x_0$};
  \fi
  \ifnum\i=1
    \node[below] at (\x,\offset-0.1) {\footnotesize $\hat x_0$};
  \fi
  \ifnum\i=2
    \node[below] at (\x,\offset-0.1) {\footnotesize $x_1$};
  \fi
  \ifnum\i=3
    \node[below] at (\x,\offset-0.1) {\footnotesize $\hat x_1$};
  \fi
  \ifnum\i=4
    \node[below] at (\x,\offset-0.1) {\footnotesize $x_2$};
  \fi
}

% Connect aligned points
\foreach \i in {0,2,4} {
  \pgfmathsetmacro\x{(\i/\nfast)*\L}
  \draw[dashed,gray] (\x,0) -- (\x,\offset);
}

% Draw small curved arrows ("bows") for fast updates
\pgfmathsetmacro\xA{(0/\nfast)*\L}
\pgfmathsetmacro\xB{(1/\nfast)*\L}
\draw[->,istblue,thick,bend left=60] (\xA,\offset) to (\xB,\offset);

\pgfmathsetmacro\xC{(1/\nfast)*\L}
\pgfmathsetmacro\xD{(2/\nfast)*\L}
\draw[->,istblue,thick,bend left=60] (\xC,\offset) to (\xD,\offset);

\end{tikzpicture}
		\caption{Timescale visualization: The upper timescale represents changes in the optimization problem; the lower one represents the number of algorithm iterations until new problem data arrives.}
		\label{fig:multistep:a}
	\end{subfigure}
	\hfill
	\begin{subfigure}[t]{0.47\textwidth}
		\centering
		\begin{tikzpicture}[
  block/.style = {draw, minimum width=0.5cm, minimum height=0.5cm, align=center},
  ->, >=Latex
]

% Nodes
\node[block] (Grad1) {\small $\nabla f_k$};
\node[block, below=0.05cm of Grad1] (Grad2) {\small $\nabla f_k$};
\node[block, below=0.3cm of Grad2] (Algo) {\footnotesize $\left[ \begin{array}{c|cc} \bar A^2 & \bar A\bar B & \bar B  \\ \hline \bar C & 0 & 0 \\ \bar C\bar A & \bar C\bar A\bar B & 0 \end{array} \right]$};

% Arrows
\draw[->] (Grad1.east) -- ++(2,0) |-([yshift=-0.5em]Algo.east) ;
\draw[->] (Grad2.east) -- ++(1.5,0) |-  ([yshift=0.5em]Algo.east);
\draw[<-] (Grad1.west) --  ++(-2,0) |- ([yshift=-0.5em]Algo.west);
\draw[<-] (Grad2.west) --  ++(-1.5,0) |-  ([yshift=0.5em]Algo.west);

\node[draw=black, dashed,  inner xsep=10pt, inner ysep=5pt, fit=(Grad2)(Grad1)] (PhiBox) {};
\node[anchor=west, yshift=16] at (PhiBox.east) {$\varphi_k$};

\end{tikzpicture}
		\caption{System-theoretic visualization: multiple iterations can be represented by lifting the system representation with the number of steps performed. The oracle becomes a repeated nonlinearity.}
		\label{fig:multistep:b}
	\end{subfigure}
    \vspace{-1em}
	\caption{Illustration of a two-step method.}
	\label{fig:multistep}
	\vspace{-3mm}
\end{figure}

	We can construct general multi-step algorithms by lifting the oracle and system description. 
	Consider a first-order one-step algorithm that is expressed as
	\begin{equation}\label{eq:onestep_algo}
		\begin{aligned}
			\xi_{k+1} &= \bar{A}(\theta_k) \xi_k + \bar{B}(\theta_k) \bar{u}_k, \\
			\bar y_k &= \bar{C}(\theta_k) \xi_k,
		\end{aligned}
		\qquad 
		\bar u_k = \nabla f_k(\bar y_k), 
		\qquad
		x_k = C_x \xi_k.
	\end{equation}
	Then its corresponding $p$-step version has the realization
	\begin{subequations}\label{eq:lifted_system}
	\begin{align}
		\setlength{\arraycolsep}{3pt}
		\begin{aligned}
			\xi_{k+1} &= \underbrace{\quad \bar A^{p} \quad}_{A(\theta_k)}\; \xi_k
			+ \underbrace{\begin{bmatrix} \bar A^{p-1}\bar B \,\, & \,\, \cdots \,\, & \,\, \bar A\bar B \,\, & \,\, \bar B \end{bmatrix}}_{B(\theta_k)} \, u_k, \\[-0.5ex]
			y_k &= \underbrace{\!\begin{bmatrix} \bar C \\ \bar C\bar A \\ \vdots \\ \bar C\bar A^{p-1} \end{bmatrix}\!}_{C(\theta_k)} \xi_k
			\! + \! \underbrace{\begin{bmatrix}
					0 & \cdots & & & 0 \\[-0.5ex]
					\bar C\bar B & 0 & & & \vdots \\
					\vdots & \ddots & \ddots & & \vdots \\
					\bar C\bar A^{p-2}\bar B & \cdots & \bar C\bar B & & 0
			\end{bmatrix}}_{D(\theta_k)} u_k,
		\end{aligned}
	\end{align}
	where we omitted the explicit dependence of the one-step matrices on $\theta_k$, and where
	\begin{equation}
		y_k = \begin{bmatrix} y^1_k \\ \vdots \\ y^p_k \end{bmatrix}, \qquad
		u_k = \begin{bmatrix} \nabla f_k(y^1_k) \\\vdots \\ \nabla f_k(y^p_k) \end{bmatrix}, \qquad x_k = C_x \xi_k.
	\end{equation}
	\end{subequations}

We emphasize that multi-step formulations are central elements in time-varying optimization algorithms, and their representation within a framework such as \cref{eq:algorithm_lpv} has not yet been considered and is novel. From a control perspective, we highlight that the oracle becomes a repeated nonlinearity, see Fig.~\ref{fig:multistep:b} for an illustration, which is well-studied in robust control and the IQC literature \cite{fetzer}. This connection may be studied in future work.

The previous examples shall emphasize the generality of \cref{eq:algorithm_lpv}. With the consideration of normal cone or projection operators, also algorithms for constrained optimization can be covered. We refer to e.g. \cite{lessardDissipativity,upadhyaya,jakob_oco}. However, for ease of exposition, we focus in the following on algorithms that either have the form \cref{eq:onestep_algo} or \cref{eq:lifted_system}, i.e., one-step accelerated gradient methods and their multi-step versions. 

\subsection*{Fixed-point property}

To lay the foundation for analyzing convergence and tracking of the algorithm, it is essential to characterize the structural properties for which \cref{eq:lifted_system} admits a fixed-point $(\xi_k^\star, y_k^\star, u_k^\star)$ for which
\begin{align}\label{eq:fixpoint}
\begin{aligned}
	\xi_k^\star &= A(\theta_k) \xi_k^\star + B(\theta_k) u_k^\star \\
	y_k^\star &= C(\theta_k) \xi_k^\star + D(\theta_k) u_k^\star,
\end{aligned}
\end{align}
and where $x^\star_k = C_x \xi^\star_k$ satisfies first-order optimality $\nabla f_k(x^\star_k) = 0$. To this end, we impose the following structural conditions on the algorithms that can be considered with our analysis approach.
\begin{assum}\label{assum:fixed_point_assum}
	For all $\theta \in \Theta$, there exists a matrix $U(\theta) \in \mathbb{R}^{n_\xi \times d}$ such that $\bar{A}(\theta)$ and $\bar{C}(\theta)$ satisfy
	\begin{equation}
		\label{eq:assum:fixed_point:michalowski}
		\begin{bmatrix}
			I - \bar A(\theta) \\  \bar C(\theta) \\ C_x
		\end{bmatrix}  U(\theta) = \begin{bmatrix}
			0 \\ I_d \\ I_d
		\end{bmatrix}.
	\end{equation}
\end{assum}
\begin{prop}
	\label{assum:fixed_point_p1}
	Let a one-step algorithm \cref{eq:onestep_algo} be given and assume that for all
	$\theta \in \Theta$ the pair $(\bar A(\theta), \bar C(\theta))$ is observable. Then,
	\cref{assum:fixed_point_assum} holds if and only if for every
	fixed $\theta \in \Theta$ and every $f \in \Smlt$, \cref{eq:onestep_algo} with fixed parameter $\theta_k = \theta$ has a
	fixed-point $\xi^\star$ with $C_x \xi^\star = \bar C(\theta)\xi^\star = x^\star$, where $x^\star = \arg\min_{x} f(x, \theta)$. In that case $\bar u^\star = 0$, $\bar y^\star = x^\star$ and $\xi^\star = U(\theta) x^\star$.
\end{prop}
\begin{proof} 
	For
	a fixed parameter $\theta$, i.e.,
	constant matrices $\bar A, \bar C, U$ and a function class fixed in $\theta$, \cref{assum:fixed_point_p1} is exactly \cite[Theorem~2]{michalowski}, whose proof is given in \cite[Appendix~A.2.2]{michalowski}. 
	The arguments remain valid for all possible fixed $\theta \in \Theta$, i.e., if the recursion \eqref{eq:onestep_algo} has such an equilibrium for every fixed $\theta_k = \theta \in \Theta$, and if \cref{eq:assum:fixed_point:michalowski} holds for all $\theta\in\Theta$. 
	For such a fixed-point, $\xi^\star = \bar A(\theta)\xi^\star$ and $C_x \xi^\star = \bar C(\theta)\xi^\star = x^\star$ yield $\mathcal{O}(\theta)\xi^\star = 1_{n_\xi}\otimes x^\star$,
	where $\mathcal{O}(\theta) = \begin{bmatrix}
	\bar C(\theta)^\top & (\bar C(\theta) \bar A(\theta))^\top & \dots & (\bar C(\theta) \bar A(\theta)^{n_\xi - 1})^\top
	\end{bmatrix}^\top$ is the observability matrix. Since $\Smlt$ is translation invariant at fixed $\theta$, $x^\star$ ranges over $\R^d$ as $f$ ranges over $\Smlt$, and we have $\mathrm{ran}(1_{n_\xi}\otimes I_d) \subseteq \mathrm{ran}(\mathcal{O}(\theta))$. By observability, $\mathcal{O}(\theta)$ has constant full column rank and thus, $\xi^\star = U(\theta)x^\star$ with $U(\theta) = \mathcal{O}(\theta)^\dagger(1_{n_\xi}\otimes I_d)$ is the unique solution.
\end{proof}
For every instance of the parameter sequence $\{ \theta_k \}_{k\geq 0}$ we obtain a fixed-point triple $\{(\xi_k^\star, \bar y_k^\star, \bar u_k^\star)\}_{k\geq 0} \! =\!  \{(\xi_k^\star, x_k^\star, 0)\}_{k\geq 0}$, and thus, a moving target for the algorithm to track. Proposition~\ref{assum:fixed_point_p1} is the (one-step) parameter-varying counterpart of the well-known structural assumptions in the static setting \cite{michalowski,scherer_ebenbauer}. The multi-step case is more delicate. An equilibrium of the lifting~\eqref{eq:lifted_system} is not necessarily unique as it may result from a period-$p$ orbit of the corresponding first-order method. We therefore restrict attention to multi-step algorithms built from convergent one-step maps.
\begin{assum}\label{assum:contraction}
	For every fixed $\theta\in\Theta$ and every $f \in \Smlt$, the recursion
	$\xi_{i+1} = \bar A(\theta)\xi_i + \bar B(\theta)\nabla f(\bar C(\theta)\xi_i,\theta)$ has a globally asymptotically stable fixed-point.
\end{assum}
Whenever \cref{assum:fixed_point_assum} holds, at this fixed-point the algorithm readout coincides with the minimizer of the corresponding static problem. Note that \cref{assum:contraction} requires convergence of the one-step algorithm when applied to the static problem with a frozen $\theta$. For classical algorithms such as GD, NM, or TM, tuning rules for $\alpha,\beta,\gamma$ in terms of $m(\theta)$ and $L(\theta)$ that ensure \cref{assum:contraction} are readily available, cf.\ \cite[Table~1]{michalowski}. General algorithms can be certified with known LMI approaches from the literature on static optimization, e.g. \cite[Theorem~4]{lessard}.
\begin{prop}\label{prop:fixed_point_pstep}
	Let a one-step algorithm \cref{eq:onestep_algo} be given. Assume that for all $\theta\in\Theta$ the pair $(\bar A(\theta), \bar C(\theta))$ is observable and that \cref{assum:contraction} holds.
	Then Assumption~\ref{assum:fixed_point_assum} holds if and only if for every
	fixed $\theta \in \Theta$ and every $f \in \Smlt$, \eqref{eq:lifted_system} with frozen parameter $\theta_k=\theta$ has a
	fixed-point $\xi^\star$ with $y^\star = 1_p \otimes x^\star$ and $C_x \xi^\star=x^\star$, where $x^\star = \arg\min_{x} f(x, \theta)$. In that case $u^\star = 0$, and $\xi^\star = U(\theta) x^\star$.
\end{prop}
The proof is provided in~\cref{sec:proofs}.
Some features are worth noting. First, if a one-step method satisfies the structural conditions \eqref{eq:assum:fixed_point:michalowski}, so does each of its multi-step versions with the same $U$ and same optimal state $\xi^\star$. Moreover, whenever observability and \cref{eq:assum:fixed_point:michalowski} holds, the optimal output is uniquely determined as $y_k^\star = 1_p \otimes x_k^\star$, expressing the consensus among the intermediate steps of the multi-step method.

\begin{rem}\label{rem:constant_U}
	In several cases, \cref{assum:fixed_point_assum} is satisfied with a constant $U$ even though $\bar A(\theta), \bar C(\theta)$ are parameter-dependent. For example, if $\bar A(\theta), \bar C(\theta)$ are chosen as in the realization of \cref{exmp:accelerated_algo}, \cref{eq:assum:fixed_point:michalowski} is satisfied with $U = \left[\begin{smallmatrix} I_d \\ I_d \end{smallmatrix}\right]$ independently of $\alpha,\beta,\gamma$. That is, this choice of $U$ satisfies \cref{assum:fixed_point_assum} for all algorithms in the family mentioned therein, e.g. GD, NM, HB, and TM. Note that for the special case $\beta=\gamma=0$, which recovers GD, ($\bar A$, $\bar{C}$) is unobservable. To obtain an observable realization we can simply take $(\bar{A},\, \bar{B}(\theta_k),\, \bar{C},\, C_x) = (I_d,\, -\alpha_k I_d,\, I_d,\, I_d)$, which satisfies \cref{eq:assum:fixed_point:michalowski} with $U = I_d$.
\end{rem}

\section{IQCs for Time-Varying Gradients}\label{sec:iqcs}

IQCs have proven to be a powerful tool to describe nonlinear input-output behaviors \cite{megretzki}.  
Informally, an operator $\varphi$ satisfies an IQC defined by some linear dynamical system $\Psi$ and a symmetric matrix $M$ if for all square-summable sequences $\{y_k\}_{k\geq 0}$ and $u_k = \varphi(y_k)$, it holds that
\begin{equation*}
\sum_{k=0}^N z_k^\top M z_k \geq 0, \qquad z_k = \Psi \begin{bmatrix} y_k \\ u_k \end{bmatrix},
\end{equation*}
for all $N \geq 0$. For an in-depth discussion of this topic we refer to the extensive IQC literature, e.g. \cite{veenmann,boczar} and all references therein. Conceptually, we will use this concept as a surrogate description of the nonlinear relation $u_k = \varphi_k(y_k)$. Many IQCs have been derived and employed for gradients of strongly convex and smooth functions \cite{lessard,scherer_ebenbauer,colombino2018}. This section applies this concept also to gradients of parameter-varying functions in $\Smlt$. We leverage filters $\Psi$ that are LPV systems themselves, which connects to the notion of parameter varying IQCs (cf. \cite{pfifer} for a discussion of LPV IQCs).

\subsection{Pointwise IQCs}\label{sec:iqc:pointwise}
As the name suggests, pointwise IQCs are constraints that hold pointwise in time. 
\begin{prop}[Sector IQC]\label{prop:sector_iqc}
    Let $f \in \Smlt$. Let $x_k$ be arbitrary for some $k\geq 0$ and let $x_k^\star$ be the solution to $\nabla f_k(x_k^\star)=0$. Define
	\begin{subequations}\label{eq:prop:sectorIQC}
	\begin{equation}
		z_k = \Psi_\theta
		\begin{bmatrix}
            x_k - x_k^\star \\ \nabla f_k(x_k)
		\end{bmatrix}, \quad  
		\Psi_\theta = \left[ \begin{array}{c|cc}
			0 & 0 & 0 \\ \hline 0 & L(\theta) I_d & -I_d \\ 0 & -m(\theta) I_d & I_d
		\end{array}
		\right].
	\end{equation}
	Then it holds
	\begin{equation}
		z_k^\top M z_k \geq 0, \qquad M = \begin{bmatrix}
			0 & I_d \\
			I_d & 0 \\
		\end{bmatrix}.
	\end{equation}
	\end{subequations}
\end{prop}

\cref{prop:sector_iqc} establishes a sector constraint on the gradient $\nabla f_k(\cdot)$ when it is held fixed in time. It is the straightforward extension of \cite[Lemma 6]{lessard} for varying convexity and Lipschitz constants, which follows naturally from the fact that the time-variation of $f$ does not need to be considered.
Pointwise IQCs are simple to handle in convergence proofs, as we will show in the next section, however, it is well known that \cref{eq:prop:sectorIQC} characterizes all nonlinearities that are sector bounded, while gradients of convex functions are in turn known to be slope restricted, which is a stronger property~\cite{lessardDissipativity}. This over-approximation has been shown to be a source of conservatism. We therefore aim to develop a more expressive IQC description, for which we need to consider dynamic filters.

\subsection{Variational IQCs}\label{sec:iqc:vIQC}

Due to the time-dependence of the stationary point $x_k^\star$ and the function landscape $f(\cdot,\theta_k)$, available results on dynamic IQCs for convex functions are not applicable to the time-varying case. Therefore, we develop an extension of the IQC concept with the consideration of the variational measures \cref{eq:variational}, for which we introduce the terminology \emph{variational IQCs}. Define 
\begin{equation}\label{eq:f_hat}
    \hat{f}_k(x) := f(x, \theta_k) - f(x_k^\star, \theta_k),
\end{equation}
so that $\hat f_k(x_k^\star)=0$. Then the following Proposition holds.

\begin{prop}[Variational IQC]\label{prop:off_by_1_iqc}
    Let $f \in \Smlt$ and $\rho \in (0,1)$ be some exponential discount factor. For some arbitrary square summable sequence $\{ x_k \}_{k\geq 0}$ and $\{ x_k^\star \}_{k\geq 0}$ being the optimal sequence solving $\nabla f_k(x_k^\star) = 0, \,\forall k\geq 0$, define the filtered sequence $\{z_k \}_{k\geq 0}$ as
	\vspace{-1mm}
	\begin{subequations}\label{eq:prop:off_by_1_iqc}
	\begin{equation}
		\label{eq:prop:off_by_1_iqc:inputs}
		z_k = \Psi_\theta
		\begin{bmatrix}
            x_k - x_k^\star \\ \nabla f_k(x_k) \\ \Delta x_k^\star \\ \Delta g_k(x_k)
		\end{bmatrix},
	\end{equation}
	where $\Delta x_k^\star, \Delta g_k(x_k)$ are measures of temporal variation as defined in \cref{eq:variational}, and $\Psi_\theta$ is an LPV system with the realization
	\begin{equation}
		\label{eq:prop:off_by_1_iqc:filter}
		\Psi_\theta =
		\left[ 
		\begin{array}{cccc|c:c:c:c} 
			0 & 0 & 0 & 0 & I_d & 0 & I_d & 0 \\ 
			0 & 0 & 0 & 0 & 0 & I_d & 0 & -I_d \\ 
			0 & 0 & 0 & 0 & -m(\theta)I_d & I_d & 0 & 0 \\ 
			0 & 0 & 0 & 0 & a(\theta)I_d & 0 & 0 & 0 \\ 
			\hline 
			-\rho^2 L(\theta)I_d & \rho^2 I_d & 0 & 0 & L(\theta)I_d & -I_d & 0 & 0\\ 
			0 & 0 & 0 & 0 & -m(\theta)I_d & I_d & 0 & 0 \\ 
			\hdashline
			0 & 0 & 0 & \rho I_d & 0 & 0 & 0 & 0 \\ 
			a(\theta) \rho I_d & 0 & 0 & 0 & 0 & 0 & 0 & 0 \\ 
			\hdashline
			0 & 0 & \rho I_d & 0 & 0 & 0 & 0 & 0 \\ 
			-m(\theta)\rho I_d & \rho I_d & 0 & 0 & 0 & 0 & 0 & 0 
		\end{array} 
		\right]_{\zeta_0=0}
	\end{equation}
	where $a(\theta) := \sqrt{\frac{m(\theta)(L(\theta)-m(\theta))}{2}}$. Then, for all $N\geq 0$ and $(m_k, L_k) := (m(\theta_k), L(\theta_k))$, it holds
	\vspace{-1mm}
	\setlength{\jot}{-2pt}
	\begin{align}
		\begin{aligned}\label{eq:prop:off_by_1_iqc:ineq}
		&\sum_{k=0}^N \rho^{-2k} z_k^\top \hat{M} z_k \\[-1ex]
		&\quad \quad \geq \sum_{k=1}^N \rho^{-2(k-1)} \left( (L_{k-1} - m_{k-1})\hat{f}_{k-1}(x_{k-1}) - (L_{k} - m_{k})\hat{f}_{k}(x_{k-1}) \right)
		\end{aligned}
	\end{align}
	\setlength{\jot}{3pt}
	with
	\vspace{-1mm}
	\begin{equation}\label{eq:prop:off_by_1_iqc:M}
		\hat{M} = \mathrm{blkdiag}\left( \frac{1}{2}\begin{bmatrix} 0 & I_d \\ I_d & 0 \end{bmatrix}, \begin{bmatrix} I_d & 0 \\ 0 & -I_d \end{bmatrix}, \frac{1}{2}\begin{bmatrix} I_d & 0 \\ 0 & -I_d \end{bmatrix} \right).
	\end{equation}
	\end{subequations}
\end{prop}
The proof is provided in \cref{sec:proofs}.
\cref{prop:off_by_1_iqc} highlights a fundamental departure from the classical IQC notion. First, \cref{eq:prop:off_by_1_iqc:inputs} incorporates explicit terms for the change in the minimizer $\Delta x_k^\star$ and the variation in the gradient $\Delta g_k(x_k)$ as input to the filter; second, since $\Delta x_k^\star$, $\Delta g_k(x_k)$ are generally not square summable, \cref{eq:prop:off_by_1_iqc:ineq} constrains the integral quadratic term to be larger than the function variation, also relating to the temporal change, rather than zero. It does therefore not present an IQC in the classical sense and hence, we introduce the phrase variational IQC. The conceptual difference with \cref{prop:sector_iqc} is illustrated in Fig.~\ref{fig:iqc_comparison}. Note that despite $\Delta x_k^\star$ and $\Delta g_k(x_k)$ being inputs to the filter, we will not require knowledge of their realizations; instead, we use this formulation as a surrogate to establish our bounds in the next section.
\begin{figure}[ht]
	\centering
	\begin{subfigure}[t]{0.47\textwidth}
		\centering
		\begin{tikzpicture}[
  block/.style = {draw, minimum width=0.5cm, minimum height=0.5cm, align=center},
  split/.style={draw, fill, circle, minimum size=0.25em, inner sep=0pt},
  ->, >=Latex
]

% Nodes
\node[block] (Grad) { $\nabla f_k$};
\node[block] (Psi) at ($(Grad.north) + (3.5em,1.5em)$) {$\Psi_\theta$};

\node[split] (splitX) at ($(Grad.west) + (-1em,0)$) {};
\node[split] (splitU) at ($(Grad.east) + (0.5em,0)$) {};

\draw[->] ($(Grad.west) + (-2.5em,0)$) -- node[left, pos=0.05] { $y_k$} ($(Grad.west)$);
\draw[->] ($(Grad.east)$) -- node[right, pos=0.95] { $u_k$} ($(Grad.east) + (3em,0)$);

\draw[->] (splitX) |- ($(Psi.west)+(0em,0.2em)$);
\draw[->] (splitU) |- ($(Psi.west)+(0em,-0.2em)$);

\draw[->] (Psi.east) -- ++(1.5em,0) node[right] { $z_k$};

\end{tikzpicture}
		\caption{Pointwise IQC: $z_k$ is realized as the filtered signal of the input and output of the gradient.}
		\label{fig:iqc_comparison:sector}
	\end{subfigure}
	\hfill
	\begin{subfigure}[t]{0.47\textwidth}
		\centering
		\begin{tikzpicture}[
  block/.style = {draw, minimum width=0.5cm, minimum height=0.5cm, align=center},
  split/.style={draw, fill, circle, minimum size=0.25em, inner sep=0pt},
  ->, >=Latex
]

% Nodes
\node[block] (Grad) { $\nabla f_k$};

% Move Ψ further to the right by +2em
\node[block, minimum height=1cm, minimum width=0.8cm] (Psi) 
    at ($(Grad.north) + (7em,1.5em)$) { $\Psi_\theta$};

\node[split] (splitX) at ($(Grad.west) + (-1em,0)$) {};
\node[split] (splitU) at ($(Grad.east) + (0.7em,0)$) {};

% ---------------------------------------------
% Define Psi input coordinates (in1 = lower, in2 = slightly above)
% ---------------------------------------------
\coordinate (PsiIn1) at ($(Psi.west)+(0em,-1.2em)$);
\coordinate (PsiIn2) at ($(Psi.west)+(0em, -0.2em)$);

% ---------------------------------------------
% New splitU2 aligned with Ψ input 1
% ---------------------------------------------
\coordinate (splitU2pos) at ($(splitU)!(PsiIn1)!(splitU)$);
\node[split] (splitU2) at (splitU2pos) {};

% ---------------------------------------------
% Δg_k aligned with Ψ input 2
% ---------------------------------------------
\coordinate (DeltaBase) at ($(splitU2)!(PsiIn2)!(splitU2)$);
\coordinate (DeltaPos) at ($(DeltaBase)+(2em,0)$);

\node[block] (Delta) at (DeltaPos) {\footnotesize $\Delta g_k$};

% ---------------------------------------------
% Arrows
% ---------------------------------------------
\draw[->] ($(Grad.west) + (-2em,0)$) -- node[left, pos=0.05] { $y_k$} ($(Grad.west)$);
\draw[->] ($(Grad.east)$) -- node[right, pos=0.95] {$u_k$} ($(Grad.east) + (3em,0)$);

% splitX → upper Ψ input
\draw[->] (splitX) |- ($(Psi.west)+(0em,0.7em)$);

% splitU → splitU2
\draw[-] (splitU) -- (splitU2);

% splitU2 → Ψ input 1
\draw[->] (splitU2) |- (PsiIn1);

% splitU2 → Δg_k
\draw[->] (splitU2) |- (Delta.west);

% Δg_k → Ψ input 2
\draw[->] (Delta.east) -- (PsiIn2);

% Ψ output
\draw[->] (Psi.east) -- ++(1.5em,0) node[right] { $z_k$};

%  input (unchanged)
\draw[->] ($(Psi.west)+(-3.1cm,1.2em)$) -- 
    node[left, yshift=0.2em ,pos=0.0] {\small \textcolor{black}{$\Delta x_k^\star$}}
    ($(Psi.west)+(0em,1.2em)$);

\end{tikzpicture}
		\caption{Variational IQC: $z_k$ is realized as the filtered signal of the input and output of the gradient, and the exogenously caused variational signals $\Delta x_k^\star$ and $\Delta g_k$.}
		\label{fig:iqc_comparison:vIQC}
	\end{subfigure}
	\vspace{-1em}
	\caption{Illustration of the difference between the pointwise and variational IQC filters.}
	\label{fig:iqc_comparison}
	\vspace{-1em}
\end{figure}

\begin{rem}
	The weighting factor $\rho^{-2k}$ serves to establish linear convergence results when using the IQC for algorithm analysis, which is a common strategy in exponential stability analysis; we refer to the literature on $\rho$-hard IQCs \cite{lessard,boczar,scherer_ZF}.
\end{rem}

\begin{rem}
	The right hand side of \cref{eq:prop:off_by_1_iqc:ineq} can be interpreted as a change in the objective function landscape. Note that in the special case of constant convexity and smoothness parameters $m, L$, respectively, and assuming differentiability of $f$ w.r.t. $\theta$, a first-order Taylor approximation would yield
	\begin{equation*}
		f(x_{k-1}, \theta_k) - f(x_{k-1}, \theta_{k-1}) \approx \nabla_\theta f(x_{k-1},\theta_{k-1})^\top (\theta_k - \theta_{k-1}).
	\end{equation*}
	This term corresponds to what in literature is often referred to as the temporal variability of the optimization problem \cite{bianchin_internalmodel}. 
\end{rem}

\begin{rem}
	Note that for static optimization problems, i.e. where $x_k^\star$, $f(\cdot, \theta_k)$, and $\nabla f_k$ become time-invariant, \cref{eq:prop:off_by_1_iqc} recovers the weighted off-by-one IQC from \cite[Lemma 8]{lessard}. We note as well that one might consider formulating ``variational off-by-$j$'' IQCs to pave the way to a generalization of Zames-Falb IQCs \cite[Lemma 8]{lessard}. This would include the consideration of variational measures of higher-order, i.e. variations across multiple timesteps. We leave this development to future works and focus on \cref{eq:prop:off_by_1_iqc} to develop our tracking bounds. Note that in \cite{scherer_ebenbauer} it was proven that off-by-one IQCs are already sufficient to certify the best possible convergence rate in the static setting, obtained by the triple momentum algorithm, and thus, their use has become standard in the IQC algorithm analysis literature.
\end{rem}

\section{Tracking Bounds for First-Order Algorithms}\label{sec:tracking}

In the following, we present a systematic approach to establish bounds on the tracking error. Instead of $\| x_k - x_k^\star \|$, we focus on deriving bounds on the distance between the algorithmic state and the optimal trajectory, i.e., $\| \xi_k - \xi_k^\star \|$. We will therefore also define a variational measure in terms of the optimal state
\begin{equation}
 \Delta \xi_k^\star = \xi_k^\star - \xi_{k+1}^\star.
\end{equation}
Throughout, we assume that $A(\theta),B(\theta),C(\theta),D(\theta)$ carry the structure of \cref{eq:lifted_system}, reducing to $(A(\theta),B(\theta),C(\theta),D(\theta))=(\bar{A}(\theta),\bar{B}(\theta),\bar{C}(\theta),0)$ for $p=1$. We further assume that $(\bar{A}(\theta),\bar{C}(\theta))$ is observable for all $\theta\in\Theta$ and that \cref{assum:fixed_point_assum} holds, together with \cref{assum:contraction} whenever $p>1$.

\subsection{Tracking Bounds with Pointwise IQCs}\label{sec:tracking:pIQC}

Given a $p$-step algorithm~\cref{eq:lifted_system}, we apply \cref{prop:sector_iqc} to every component in the input-output pair $(y_k^i, u_k^i)$, $i \in \mathbb{I}_p$, with $u_k^i = \nabla f_k(y_k^i)$. Define
\begin{equation*}
    z^i_k = \left[ \begin{array}{c|cc}
			0 & 0 & 0 \\ \hline 0 & L(\theta_k) I_d & -I_d \\ 0 & -m(\theta_k) I_d & I_d
		\end{array}
		\right]
        \begin{bmatrix} y^i_k-x_k^\star \\ u^i_k \end{bmatrix},
\end{equation*}
then by \cref{prop:sector_iqc}, we have $(z^i_k)^\top M z^i_k \geq 0$ for all $k\geq 0$ and each $i\in\mathbb{I}_p$. Let us collect all $z^i_k$ and stack them into one big vector $z_k = \mathrm{vec}(z_k^1, \dots, z_k^p)$, and recall that $y_k^\star = (1_p \otimes I_d)x_k^\star$. Then, we can define suitable matrices $D_\Psi^y$ and $D_\Psi^u$ with which we can write down the more compact relation
\begin{align}\label{eq:unaugmented_output}
 z_k = 
 \begin{bmatrix}
    D_\Psi^y(\theta_k) & D_\Psi^u(\theta_k)
 \end{bmatrix}
 \begin{bmatrix}
    y_k - y_k^\star \\ u_k
 \end{bmatrix},
\end{align}
where the block rows of $D_\Psi^y$ and $D_\Psi^u$ select the respective components of $y_k - y_k^\star$ and $u_k$ to realize each $z_k^i$.

We now use that $y_k - y_k^\star = C(\theta_k) (\xi_k - \xi_k^\star) + D(\theta_k) u_k$ to obtain
\begin{align}\label{eq:augmented_output}
    z_k &= \underbrace{\begin{bmatrix}
            D_\Psi^y(\theta_k) C(\theta_k) & D_\Psi^y(\theta_k) D(\theta_k) + D_\Psi^u(\theta_k)
           \end{bmatrix}}_{=:  \begin{bmatrix} \hat{C}(\theta_k) & \hat{D}(\theta_k) \end{bmatrix}}
           \begin{bmatrix}
            \xi_k - \xi_k^\star \\ u_k
           \end{bmatrix}.
\end{align}
We can interpret \cref{eq:augmented_output} as the augmented output of algorithm \cref{eq:system_lpv}, which neglects the oracle $\varphi$ and instead becomes an open-loop system with input $u_k$ and output $z_k$, see Fig.~\ref{fig:algorithm_pIQC}. We can then state the following theorem.

\begin{figure}
\centering
\begin{tikzpicture}[
  >=Latex,
  block/.style = {
    draw, minimum width=0.8cm, minimum height=0.6cm,
    align=center, fill=white
  },
  gradblock/.style = {
    draw=gray!70,
    text=gray!70,
    minimum width=0.9cm,
    minimum height=0.6cm,
    align=center,
    fill=white
},
  stacklayer/.style = {
    draw=gray!70,
    fill=white,
    minimum width=0.9cm, minimum height=0.6cm
},
stacklayertwo/.style = {
    draw=black,
    fill=white,
    minimum width=0.8cm, minimum height=0.6cm},
  % original arrow styles
  vecarrow/.style = {->, double, line width=0.3pt, >=Latex, double distance=0.8pt},
  vecarrowop/.style = {<-, double, line width=0.3pt, >=Latex, double distance=0.8pt},
  % gray internal arrow (with tip)
  vecarrowgray/.style = {
    ->, 
    double,
    line width=0.3pt,
    double distance=0.8pt,
    >=Latex,
    gray!70    % <--- change this
},
doublegrayplain/.style = {
    double,
    line width=0.3pt,
    double distance=0.8pt,
    gray!70    % <--- change this
},
  split/.style = {circle, fill=black, minimum size=0.25em, inner sep=0pt}
]

% Algo block
\node[block] (Algo) at (0,0) {
  $\begin{aligned}
    \xi_{k+1} &= A(\theta_k)\,\xi_k + B(\theta_k)\,u_k \\
    y_k &= C(\theta_k)\,\xi_k + D(\theta_k)\,u_k
  \end{aligned}$
};

% ===== Restored GRADIENT STACK =====
\begin{scope}[shift={($(Algo.north)+(0,2em)$)}]
  % back layers (light gray)
  \foreach \y in {0.12, 0.06} {
    \node[stacklayer] at (\y,\y) {};
  }
  % main gradient block
  \node[gradblock] (Grad) at (0,0) {$\nabla f_k$};
\end{scope}

% Split nodes
\node[split] (splitL) at ($(Grad.west)+(-1.6em,0)$) {};
\node[split] (splitR) at ($(Grad.east)+(1.6em,0)$) {};

% ==== RIGHT SIDE: y_k path (Algo -> split -> Grad) ====

% Algo -> splitR (DOUBLE LINE, NO TIP)
\draw[double, line width=0.3pt, double distance=0.8pt]
  ([yshift=0em]Algo.east) -- ++(1.5em,0) |- (splitR)
  node[pos=0.25, right] {$y_k$};

% splitR -> Grad (GRAY ARROW, ENTERS GRAD)
\draw[vecarrowgray] (splitR) -- (Grad.east);

% ==== LEFT SIDE: u_k path (Grad -> splitL (no tip) -> Algo) ====

% Grad -> splitL  (GRAY DOUBLE LINE, NO ARROW TIP)
\draw[doublegrayplain] (Grad.west) -- (splitL);

% Algo <- splitL : first left then down-right (NO TIP), vector-pointing INTO Algo
\draw[vecarrowop, line width=0.3pt, double distance=0.8pt]
  ([yshift=0em]Algo.west) -- ++(-1.5em,0) |- (splitL)
  node[pos=0.25, left] {$u_k$};

\begin{scope}[shift={($(Grad.north east)+(4.8em,1.8em)$)}]
  \foreach \y in {0.12,0.06} {
    \node[stacklayertwo] at (\y,\y) {};
  }
  \node[block] (Psi) at (0,0) {$\Psi_\theta$};
\end{scope}

\draw[vecarrow] (splitL)  |- ($(Psi.west)+(0,0.3em)$);
\draw[vecarrow] (splitR) |- ($(Psi.west)+(0,-0.3em)$);

\draw[vecarrow] (Psi.east) -- ++(1.6em,0) node[right] {$z_k$};

\end{tikzpicture}
\caption{Augmented algorithm. The influence of the gradient is replaced by an IQC.}
\label{fig:algorithm_pIQC}
\vspace{-3mm}
\end{figure}

\begin{theorem}
    \label{thm:tracking:sector}
    Consider algorithm \cref{eq:lifted_system}. Assume the rate bound \cref{eq:param_rate_bound} holds and define the corresponding set $\Vtnu$ as in \cref{eq:consistent_polytope}. Using $\Psi_\theta$ from \cref{eq:prop:sectorIQC}, form the augmented output matrices $\hat{C}, \hat{D}$ as in \cref{eq:augmented_output}. If there exists a continuous parameter-dependent symmetric matrix $P: \Theta \rightarrow \mathbb{S}^{n_\xi}$ for which $P(\theta) \succ 0$ for all $\theta \in \Theta$ and a continuous function $\lambda_p: \Theta \rightarrow \R_{\geq 0}^p$, such that for $\rho \in (0,1)$ it holds
    \begin{align}\label{eq:thm1:lmi}
    \begin{aligned}
        \begin{bmatrix}
        I & 0  \\
        {A}(\theta) & {B}(\theta) \\
        \hat{C}(\theta) & \hat{D}(\theta)  \\
        \end{bmatrix}^\top
        \begin{bmatrix}
        \textcolor{black}{-\rho^2 P(\theta)} & 0 & 0 \\
        0 & P(\theta_+) & 0 \\
        0 & 0 & M_\lambda(\theta) \\
        \end{bmatrix}
        \begin{bmatrix}
        I & 0  \\
        {A}(\theta) & {B}(\theta) \\
        \hat{C}(\theta) & \hat{D}(\theta)  \\
        \end{bmatrix}
        \preceq 0, \\ \text{ with } \theta_+ := \theta + \Delta \theta, \quad \forall \left( \theta, \Delta \theta \right) \in \Vtnu,
    \end{aligned}
    \end{align} 
    where $M_\lambda(\theta) = \mathrm{diag}(\lambda_p(\theta)) \otimes M$,
    then for any $\xi_0$ and $k\geq 0$ we have
    \begin{equation}
        \label{eq:thm1:bound}
        \norm{\xi_k - \xi^\star_k} \leq c \rho^k \norm{\xi_0 - \xi^\star_0} + c \sum_{t=1}^{k} \rho^{k-t} \norm{\Delta \xi^\star_{t-1}},
    \end{equation}
    with $c = \sqrt{\overline{\lambda} / \underline{\lambda}}$, $\overline{\lambda}:= \max_{\theta \in \Theta} \lambda_{\max}\left( P(\theta) \right)$, $\underline{\lambda}:= \min_{\theta \in \Theta} \lambda_{\min} \left( P(\theta) \right)$.
\end{theorem}
\begin{proof}
    To streamline the notation we define ${A}_k := {A}(\theta_k)$ and equivalently for ${B}_k, \hat{C}_k, \hat{D}_k, P_k,\lambda_{p,k}$. 
    Since $P_{k+1}\succ 0\,\, \forall k$, we can leverage the triangle inequality to obtain
    \begin{align}
        \begin{aligned}
            \norm{\xi_{k+1} - \xi_{k+1}^\star}_{P_{k+1}} &= \norm{\xi_{k+1} - \xi_{k}^\star + \xi_{k}^\star - \xi_{k+1}^\star}_{P_{k+1}}
            \\ 
            \label{eq:proof_lpv_1} 
            &\leq \norm{\xi_{k+1} - \xi_{k}^\star}_{P_{k+1}} + \norm{\xi^\star_{k+1} - \xi_{k}^\star}_{P_{k+1}}.
        \end{aligned}
    \end{align}
    Now consider the first term on the right hand side in (\ref{eq:proof_lpv_1}). We use the algorithm iteration and the fixed-point property $\xi^\star_k = A \xi_k^\star$ to establish the relation
        \begin{align}\label{eq:thm1:proof:1}
            \begin{aligned}
            \norm{\xi_{k+1} - \xi_{k}^\star}_{P_{k+1}}^2 &= \norm{{A}_k \xi_k + {B}_k u_k - \xi_{k}^\star}_{P_{k+1}}^2\\
            &= \norm{{A}_k (\xi_k - \xi_k^\star) + {B}_k u_k}_{P_{k+1}}^2.
            \end{aligned}
        \end{align}
    If we now left and right multiply \cref{eq:thm1:lmi} with $\mathrm{vec}(\xi_k-\xi_k^\star, u_k)$ and use \cref{eq:augmented_output}, we obtain
    \begin{equation}\label{eq:thm1:proof:2}
        - \rho^2 \norm{\xi_k - \xi_k^\star}_{P_k}^2 + \norm{{A}_k (\xi_k - \xi_k^\star) + {B}_k u_k}_{P_{k+1}}^2   + \sum_{i=1}^p \lambda_{p,k}^i (z_k^i)^\top M z_k^i \leq 0.
    \end{equation}
    By \cref{prop:sector_iqc}, the quadratic term is nonnegative for every $i\in\mathbb{I}_p$. By combining \cref{eq:thm1:proof:1} and \cref{eq:thm1:proof:2} and taking the square root, we conclude
    \begin{equation}
        \label{eq:proof_lpv_2}
        \norm{\xi_{k+1} - \xi_{k}^\star}_{P_{k+1}} \leq \rho \norm{\xi_k - \xi_k^\star}_{P_k}.
    \end{equation}
    Plugging (\ref{eq:proof_lpv_2}) into (\ref{eq:proof_lpv_1}), we obtain
    \begin{equation}
        \label{eq:proof_lpv_3}
        \norm{\xi_{k+1} - \xi_{k+1}^\star}_{P_{k+1}} \leq \rho \norm{\xi_k - \xi_k^\star}_{P_k} + \norm{\xi^\star_{k+1} - \xi_{k}^\star}_{P_{k+1}}.
    \end{equation}
    Finally, applying (\ref{eq:proof_lpv_3}) recursively from starting index $k$ to zero results in
    \begin{align}
    \begin{aligned}
        \label{eq:proof_lpv_4}
        \norm{\xi_k - \xi^\star_k}_{P_k} &\leq \rho^k \norm{\xi_0 - \xi^\star_0}_{P_0} + \sum_{t=0}^{k-1} \rho^{k-t-1} \norm{\xi^\star_{t+1} - \xi^\star_{t}}_{P_{t+1}} \\
        &= \rho^k \norm{\xi_0 - \xi^\star_0}_{P_0} +  \sum_{t=1}^{k} \rho^{k-t} \norm{\Delta \xi^\star_{t-1}}_{P_{t}}.
    \end{aligned}
    \end{align}
    Note that since $P$ is continuous and positive definite on the compact set $\Theta$, $\underline{\lambda}$ and $\overline{\lambda}$ exist and are positive. Hence,
    \begin{equation}
        \label{eq:condition_P}
        \underline{\lambda} \norm{x}^2 \leq x^\top P(\theta) x \leq  \overline{\lambda} \norm{x}^2, \quad \forall \theta\in \Theta, x \in \R^n,
    \end{equation}
    which we apply to (\ref{eq:proof_lpv_4}), and divide by $\sqrt{\underline{\lambda}}$
    to arrive at \cref{eq:thm1:bound}.
\end{proof}

\cref{thm:tracking:sector} is, to the best of the authors' knowledge, the first result that formulates the analysis of tracking algorithms in time-varying optimization as an LPV-IQC problem. 
We remark that \cref{eq:thm1:lmi} and the LMI condition from \cite[Theorem 4]{lessard} look similar, though they differ profoundly due to the parameter-dependence of all involved matrices, including the consideration of the rate bounds within $\Vtnu$. The static bound obtained in \cite{lessard} is recovered for $\xi^\star_k \equiv \mathrm{const}$. 

\begin{corollary}\label{corr:iterates}
    Consider algorithm~\cref{eq:lifted_system} and let \cref{assum:fixed_point_assum} hold with constant $U$. Then, if the algorithm is initialized with $\xi_0 = U x_0$, the following iterate bound follows from \cref{thm:tracking:sector}
    \begin{equation}
        \label{eq:thm1:bound:x}
        \norm{x_k - x^\star_k} \leq \bar{c} \rho^k \norm{x_0 - x^\star_0} + \bar{c} \sum_{t=1}^{k} \rho^{k-t} \norm{\Delta x^\star_{t-1}},
    \end{equation}
    with $\bar{c} = c \| C_x \|_2 \| U \|_2$, and $c$ defined in \cref{thm:tracking:sector}.
\end{corollary}
\begin{proof}
    By \cref{assum:fixed_point_p1} or \ref{prop:fixed_point_pstep} the unique optimal algorithmic state is $\xi_k^\star = U x_k^\star$. Hence, for all $k$ and all $t$ the bounds
    \begin{align}\label{eq:road_to_rm52}
        \|x_k - x_k^\star\| &= \|C_x(\xi_k - \xi_k^\star)\| \leq \|C_x\|_2 \|\xi_k - \xi_k^\star\| \\
        \| \xi_0 - \xi_0^\star \| &= \| U(x_0 - x_0^\star) \| \leq \| U \|_2 \| x_0 - x_0^\star \|\notag \\
       \| \Delta \xi_{t-1}^\star \| = \| \xi_t^\star - \xi_{t-1}^\star \| &= \| U (x_t^\star - x_{t-1}^\star) \| \leq \| U \|_2 \| x_t^\star - x_{t-1}^\star \| = \| U \|_2 \| \Delta x_{t-1}^\star \| \notag
    \end{align}
    hold. Eq. \cref{eq:thm1:bound} together with \cref{eq:road_to_rm52} implies \cref{eq:thm1:bound:x}.
\end{proof}
\cref{thm:tracking:sector} certifies the algorithm's behavior in terms of the algorithmic state, \cref{corr:iterates} translates this into a bound on the iterate $x_k$. In the special case of a constant $U$, cf. \cref{rem:constant_U}, \cref{corr:iterates} shows that we recover the particularly simple form of \cref{eq:tracking_bound} with $c_1=c_2=c\|C_x\|_2\|U\|_2$ and $\delta_t = \|\Delta x_{t-1}^\star\|$. 
If \cref{assum:fixed_point_assum} holds instead with a $\theta$-dependent $U(\theta)$, and the algorithm is initialized with $\xi_0 = U(\theta_0)x_0$, the iterate bound \cref{eq:tracking_bound} can be recovered with $c_1 = c\|C_x\|_2 \bar{U}$, $c_2 = c\|C_x\|_2$ and $\delta_t = \| U(\theta_t) x_t^\star - U(\theta_{t-1}) x_{t-1}^\star\|$, where $\bar{U} := \mathrm{sup}_{\theta\in\Theta} \| U(\theta)\|_2$ is finite\footnote{Recall that $U(\theta) = \mathcal{O}^{\dagger}(\theta)(1_{n_{\xi}} \otimes I_d)$, $\bar{A}(\theta)$ and $\bar{C}(\theta)$ are continuous in $\theta$, $\Theta$ is compact, and ($\bar{A}(\theta), \bar{C}(\theta)$) is observable for all $\theta\in\Theta$. Hence, $\mathcal{O}(\theta)$ is continuous and has constant full rank on $\Theta$, thus $\mathcal{O}^{\dagger}(\theta)$ is continuous and $\bar{U}$ is finite}.

\begin{rem}
    The parametrization of the Lyapunov matrix $P(\theta)$ is a common technique to reduce the conservatism of constant Lyapunov matrices in rate-bounded LPV stability tests \cite{deSouza_consistentPolytope,amato}. Note that the existence of a continuous matrix-valued function $P: \Theta \rightarrow \R^{n \times n}$ satisfying \cref{eq:thm1:lmi} for all $\theta$ cannot generally be inferred a priori from the problem data alone; it is therefore an explicit assumption on the decision variable. Moreover, checking its existence with the LMI is also nontrivial, since, unless $\Theta$ is polytopic and the dependence of all involved matrices on $\theta$ is affine, verifying \cref{eq:thm1:lmi} is an infinite-dimensional problem. In practice, a common approach is to linearly parametrize $P$, e.g.\ by choosing $N_P \in \mathbb{N}$ continuous ansatz functions $\phi_i: \Theta \rightarrow \R$ yielding $P(\theta) = P_0 + \sum_{i=1}^{N_P} P_i \phi_i(\theta)$ with decision variables $P_i$, and to grid over $\Vtnu$ \cite{wuPackard_boundedRateLyap,akparian, pfifer}. Note that a solution to this gridded problem does not guarantee feasibility of the original (infinite-dimensional) LMI.
\end{rem}

\subsection{Tracking Bounds with Variational IQCs}\label{sec:tracking:vIQC}

Following the discussion on conservatism of pointwise IQCs in section~\ref{sec:iqc:pointwise}, we now aim for a tracking bound that leverages the variational IQC. We proceed as before: \cref{prop:off_by_1_iqc} holds for every component in the pair $(y_k^i, u_k^i)$, $i \in \mathbb{I}_p$, with the additional filter inputs $\Delta x_k^\star$ and $\Delta g_k(y_k^i)$ (recall \cref{eq:prop:off_by_1_iqc:inputs}). Define the filtered sequences $\{z^i_k\}_{k\geq 0}$ as
\vspace{-2mm}
\begin{equation}\label{eq:psi_i:offbyone}
    z^i_k = 
         \left[ \begin{array}{c|cccc}
A_{\Psi^i} & B_{\Psi^i}^y & B_{\Psi^i}^u & B_{\Psi^i}^{\Delta x^\star} & B_{\Psi^i}^{\Delta g} \\ \hline
C_{\Psi^i} & D_{\Psi^i}^y & D_{\Psi^i}^u & D_{\Psi^i}^{\Delta x^\star} & D_{\Psi^i}^{\Delta g}
\end{array} \right]_{\zeta_0^i}
        \begin{bmatrix}
            y^i_k - x_k^\star \\ u_k^i \\ C_x  \Delta \xi_k^\star \\ \Delta g_k(y^i_k)
		\end{bmatrix},
\end{equation}
where the matrices $A_{\Psi^i}$, $B_{\Psi^i}^\star$, $C_{\Psi^i}$, $D_{\Psi^i}^\star$ are defined as in \cref{eq:prop:off_by_1_iqc:filter} with some arbitrary filter initialization $\zeta_0^i$. Note that we substituted $\Delta x_k^\star = C_x \Delta \xi_k^\star$, and the dependence on the parameters has been omitted for readability. Now collect again all $z^i_k$ to stack them into one big vector $z_k = \mathrm{vec}(z_k^1, \dots, z_k^p)$ and let us additionally define the virtual input signal

\begin{equation} 
w_k := \mathrm{vec}(\Delta g_k(y^1_k),\dots, \Delta g_k(y^p_k)),
\end{equation}
representing the gradient variation evaluated along all intermediate points of the multi-step algorithm.
Then, using straightforward system theory, we can obtain a compact filter realization by vertically stacking all filters in \cref{eq:psi_i:offbyone} and suitably permuting the input matrices, which we write compactly as
\begin{equation}\label{eq:big_IQC}
    \vspace{-2mm}
z_k = \left[ \begin{array}{c|cccc}
A_\Psi & B_\Psi^y & B_\Psi^u & B_\Psi^{\Delta x^\star} C_x & B_\Psi^{\Delta g} \\ \hline
C_\Psi & D_\Psi^y & D_\Psi^u & D_\Psi^{\Delta x^\star} C_x & D_\Psi^{\Delta g}
\end{array} \right]_{\zeta_0}
\begin{bmatrix}
    y_k - y_k^\star \\ u_k \\ \Delta \xi_k^\star \\ w_k
\end{bmatrix}.
\end{equation}
We denote the internal filter state of \cref{eq:big_IQC} by $\zeta_k$.
Here, the block rows of $B_\Psi^y$, $B_\Psi^u$, $B_\Psi^{\Delta g}$, $D_\Psi^y$, $D_\Psi^u$, and $D_\Psi^{\Delta g}$ select the respective components of $y_k - y_k^\star$, $u_k$ and $w_k$ to realize each $z_k^i$, respectively.

To proceed further, we perform a coordinate change of the algorithm from $\xi_k$ to the error coordinate $\tilde{\xi}_k : = \xi_k - \xi^\star_k$. Using the fixed-point property, we can write the mapping from $u_k$ to $\tilde{y}_k := y_k - y_k^\star$ as
\begin{align}
    \begin{aligned}
    \label{eq:generalized_algorithm_delta}
    \tilde{\xi}_{k+1} &= A(\theta_k) \tilde{\xi}_k + B(\theta_k) u_k + \Delta \xi_k^\star\\
    \tilde{y}_k &= C(\theta_k) \tilde{\xi}_k + D(\theta_k) u_k,
    \end{aligned}
\end{align}
or, equivalently,
\vspace{-2mm}
\begin{equation}\label{eq:big_plant}
\tilde y_k = \left[ \begin{array}{c|ccc}
A & B & I_d & 0 \\ \hline
C & D & 0 & 0
\end{array} \right]_{\tilde{\xi}_0}
\begin{bmatrix}
u_k \\ \Delta \xi_k^\star \\ w_k
\end{bmatrix}.
\end{equation}
Finally, use \cref{eq:big_plant} to eliminate ${y}_k - y_k^\star$ from \cref{eq:big_IQC}, which gives the augmented mapping 
\begin{align}\label{eq:augmented_plant_lpv}
    z_k = 
    \underbrace{
    \left[ \begin{array}{cc|ccc} A & 0 & B & I & 0 \\ B_\Psi^y C & A_\Psi & B_\Psi^y D + B_\Psi^u & B_\Psi^{\Delta x^\star} C_x & B_\Psi^{\Delta g} \\ \hline D_\Psi^y C & C_\Psi & D_\Psi^y D + D_\Psi^u & D_\Psi^{\Delta x^\star} C_x & D_\Psi^{\Delta g} \end{array} \right]_{\left[\begin{smallmatrix}
        \tilde{\xi}_0\\ \zeta_0
    \end{smallmatrix}\right]}}_{=:\left[ \begin{array}{c|c c c} \hat{A} & \hat{B} & \hat{B}_{\Delta \xi} & \hat{B}_{\Delta g} \\ \hline \hat{C} & \hat{D} & \hat{D}_{\Delta \xi} & \hat{D}_{\Delta g} \end{array} \right]_{\eta_0}} \begin{bmatrix}
    u_k \\ \Delta \xi_k^\star \\ w_k
    \end{bmatrix}.
\end{align}
with internal state $\eta_k := \mathrm{vec}(\tilde{\xi}_k, \zeta_k) \in \mathbb{R}^{n_\eta}$. This augmented system may be visualized by Fig.~\ref{fig:algorithm_vIQC}. We note that by \cref{eq:prop:off_by_1_iqc:filter}, the augmented plant \cref{eq:augmented_plant_lpv} is parametrized by an exponential discount factor $\rho\in(0,1)$. Finally, we are able to state the main theorem.

\begin{figure}
\centering
\begin{tikzpicture}[
  >=Latex,
  block/.style = {
    draw, minimum width=0.8cm, minimum height=0.6cm,
    align=center, fill=white
  },
  gradblock/.style = {
    draw=gray!70,
    text=gray!70,
    minimum width=0.9cm,
    minimum height=0.6cm,
    align=center,
    fill=white
},
  stacklayer/.style = {
    draw=gray!70,
    fill=white,
    minimum width=0.9cm, minimum height=0.6cm
},
stacklayertwo/.style = {
    draw=black,
    fill=white,
    minimum width=0.8cm, minimum height=1.35cm},
stacklayerthree/.style = {
    draw=black,
    fill=white,
    minimum width=0.85cm, minimum height=0.6cm
},
  % original arrow styles
  vecarrow/.style = {->, double, line width=0.3pt, >=Latex, double distance=0.8pt},
  vecarrowop/.style = {<-, double, line width=0.3pt, >=Latex, double distance=0.8pt},
  % gray internal arrow (with tip)
  vecarrowgray/.style = {
    ->, 
    double,
    line width=0.3pt,
    double distance=0.8pt,
    >=Latex,
    gray!70    % <--- change this
},
doublegrayplain/.style = {
    double,
    line width=0.3pt,
    double distance=0.8pt,
    gray!70    % <--- change this
},
doubleplain/.style = {
    -,
    double,
    line width=0.3pt,
    double distance=0.8pt
},
  split/.style = {circle, fill=black, minimum size=0.25em, inner sep=0pt}
]

% Algo block
\node[block] (Algo) at (0,0) {
  $\begin{aligned}
    \tilde \xi_{k+1} &= A(\theta_k)\tilde \xi_k + B(\theta_k)u_k + \Delta \xi_k^\star\\
    \tilde y_k &= C(\theta_k)\tilde \xi_k + D(\theta_k)u_k
  \end{aligned}$
};

% ===== Restored GRADIENT STACK =====
\begin{scope}[shift={($(Algo.north)+(0,2em)$)}]
  % back layers (light gray)
  \foreach \y in {0.12, 0.06} {
    \node[stacklayer] at (\y,\y) {};
  }
  % main gradient block
  \node[gradblock] (Grad) at (0,0) {$\nabla \tilde{f}_k$};
\end{scope}

% Split nodes for gradient
\node[split] (splitL) at ($(Grad.west)+(-1.6em,0)$) {};
\node[split] (splitR) at ($(Grad.east)+(1.6em,0)$) {};

% ==== RIGHT SIDE: y_k path (Algo -> splitR -> Grad) ====

% Algo -> splitR (DOUBLE LINE, NO TIP)
\draw[double, line width=0.3pt, double distance=0.8pt]
  ([yshift=0em]Algo.east) -- ++(1.5em,0) |- (splitR)
  node[pos=0.25, right] {$\tilde{y}_k$};

% splitR -> Grad (GRAY ARROW)
\draw[vecarrowgray] (splitR) -- (Grad.east);

% ==== LEFT SIDE: u_k path (Grad -> splitL -> Algo) ====

% Grad -> splitL  (GRAY DOUBLE LINE, NO TIP)
\draw[doublegrayplain] (Grad.west) -- (splitL);

% Algo <- splitL (NO TIP), u_k label
\draw[vecarrowop, line width=0.3pt, double distance=0.8pt]
  ([yshift=0.66em]Algo.west) -- ++(-1.5em,0) |- (splitL)
  node[pos=0.25, left] {$u_k$};

% ==== Psi block (unchanged) ====
% \begin{scope}[shift={($(Grad.north east)+(4.8em,1.8em)$)}]
%   \foreach \y in {0.12,0.06} {
%     \node[stacklayertwo] at (\y,\y) {};
%   }
%   \node[block] (Psi) at (0,0) {$\Psi_\theta$};
% \end{scope}

% --- Make Ψ taller: redefine it here ---
\begin{scope}[shift={($(Grad.north east)+(9em,1.8em)$)}]
  \foreach \y in {0.18,0.12,0.06} {  % slightly larger background stack
    \node[stacklayertwo] at (\y,\y) {};
  }
  \node[block, minimum height=1.35cm] (Psi) at (0,0) {$\Psi_\theta$};
\end{scope}

% --- Define horizontal input levels on Ψ (bottom → top) ---
\coordinate (PsiIn1) at ($(Psi.west)+(0, -0.6cm)$);  % from splitR
\coordinate (PsiIn2) at ($(Psi.west)+(0, -0.24cm)$);  % from splitL2
\coordinate (PsiIn3) at ($(Psi.west)+(0,  0.3cm)$);  % from Δg_k
\coordinate (PsiIn4) at ($(Psi.west)+(0,  0.6cm)$);  % from Δξ*_k (existing)

% ============================
% NEW splitR2
% ============================

\coordinate (splitR2) at ($(splitR)!(PsiIn1)!(splitR)$);
\node[split] (splitR2) at (splitR2) {};

% splitR → splitR2
\draw[doubleplain] (splitR) -- (splitR2);

% splitR2 → PsiIn1 (BOTTOM input)
\draw[vecarrow] (splitR2) -- (PsiIn1);

% ============================
% MOVE Δg stack to the right
% aligned with PsiIn2
% ============================

% new position: slightly right of splitR2, same height as PsiIn2
% Projection of PsiIn2 onto the vertical through splitR2
\coordinate (DeltaGbase) at ($(splitR2)!(PsiIn2)!(splitR2)$);

% Shift slightly right for visual separation
\coordinate (DeltaGpos) at ($(DeltaGbase)+(3em,0)$);

\foreach \y in {0.12, 0.06} {
    \node[stacklayerthree] at ($(DeltaGpos)+(\y,\y)$) {};
}

\node[block] (DeltaG) at (DeltaGpos) {$\Delta g_k$};

% splitR2 → Δg
\draw[vecarrow] (splitR2) |- (DeltaG.west);

% Δg → PsiIn2
\draw[vecarrow] (DeltaG.east) -- (PsiIn2);

% splitL → PsiIn3
\draw[vecarrow] (splitL) |- (PsiIn3);

% ============================================================
%       Existing Δξ*_k branch (top input)
%       (This arrow was defined elsewhere; here only attach it)
% ============================================================

% ============================================================
%       Output of Ψ
% ============================================================

\draw[vecarrow] (Psi.east) -- ++(3em,0) node[right] {$z_k$};

% ==== NEW: Extra input Δξ*_k entering Algo from the west ====

% Split placed slightly left of Algo.west
\node[split] (splitDX) at ([xshift=-3.5em, yshift=-0.66em]Algo.west) {};

% === New C_x block on same height as PsiIn4 ===
\node[block, minimum width=0.8cm, minimum height=0.6cm, rounded corners=3pt]
    (Cx) at ($(PsiIn4)+(-5cm,0)$) {$C_x$};

% Incoming arrow INTO splitDX (normal arrow, unlabeled arrow tip)
\draw[-, >=Latex]
  ([xshift=-2em]splitDX.west) -- (splitDX)
  node[pos=0, left] {$\Delta \xi^\star_k$};

% Branch 1: splitDX -> Algo
\draw[->, >=Latex]
  (splitDX) -- ([yshift=-0.66em]Algo.west);

% Branch 2: splitDX -> C_x (go up and right if necessary)
\draw[->, >=Latex]
  (splitDX) |- (Cx.west);

% Output of C_x -> Psi (horizontal, aligned with PsiIn4)
\draw[->, >=Latex]
  (Cx.east) -- (PsiIn4);

\end{tikzpicture}
\caption{Augmented interconnection with variational IQCs.}
\label{fig:algorithm_vIQC}
\vspace{-4mm}
\end{figure}

\begin{theorem}
\label{thm:tracking:variational}
Consider algorithm \cref{eq:lifted_system}. Assume the rate bound \cref{eq:param_rate_bound} holds and define the corresponding set $\Vtnu$ as in \cref{eq:consistent_polytope}. Take $\Psi_\theta$, $\hat{M}$ from \cref{eq:prop:off_by_1_iqc} for some $\rho\in(0,1)$, and form the augmented plant \cref{eq:augmented_plant_lpv}. If there exists a continuous parameter-dependent symmetric matrix $P: \Theta \rightarrow \mathbb{S}^{n_\eta}$, such that $P(\theta) \succ 0$ for all $\theta \in \Theta$, and there exists a vector $\lambda_p \in \mathbb{R}_{\geq 0}^p$ and non-negative scalars $\gamma_\xi, \gamma_g \geq 0$, such that 
    \begin{align}
    \label{eq:thm2:lmi}
    \begin{aligned}
        \begin{bmatrix}
        \star
        \end{bmatrix}^\top
        \hspace{-3pt}
        \begin{bmatrix}
        -\rho^2 P(\theta) & 0 & 0 & 0 & 0 \\
        0 & P(\theta_+) & 0 & 0 & 0 \\
        0 & 0 & M_\lambda & 0 & 0 \\
        0 & 0 & 0 & -\gamma_\xi I & 0 \\
        0 & 0 & 0 & 0 & -\gamma_g I
        \end{bmatrix}
        \hspace{-5pt}
        \begin{bmatrix}
        I & 0 & 0 & 0\\
        \hat{A} & \hat{B} & \hat{B}_{\Delta \xi} & \hat{B}_{\Delta g} \\ \hat{C} & \hat{D} & \hat{D}_{\Delta \xi} & \hat{D}_{\Delta g} \\
        0 & 0 & I & 0 \\ 0 & 0 & 0 & I
        \end{bmatrix}
        \preceq 0 & \\ \text{ with } \theta_+ := \theta + \Delta \theta, \quad \forall \left( \theta, \Delta \theta \right) \in \Vtnu,&
    \end{aligned}
\end{align}
holds,
where $M_\lambda = \mathrm{diag}(\lambda_p) \otimes \hat{M}$,
then for any $\xi_0$ and for all $k\geq 0$ there exist constants $c_1 := \overline{\lambda}/\underline{\lambda}$ and $c_2 := 1/\underline{\lambda}$ with $\overline{\lambda}, \underline{\lambda}$ as in (\ref{eq:condition_P}), such that
\vspace{-2mm}
\begin{multline}\label{eq:thm2:bound}
    \norm{\xi_k - \xi_k^\star}^2 \leq c_1 \rho^{2k} \norm{\xi_0 - \xi_0^\star}^2 + c_2 \sum_{t=1}^{k} \rho^{2(k-t)} \left( \gamma_{\xi} \norm{\Delta \xi^\star_{t-1}}^2 + \gamma_g \sum_{i=1}^{p} \| \Delta g_{t-1}(y_{t-1}^i) \|^2  \right) \\[-1ex] + c_2 \sum_{i=1}^{p} \sum_{t=1}^{k-1} \rho^{2(k-t)} \lambda_p^i \left( (L_{t} - m_{t}) \hat{f}_{t}(y^i_{t-1}) - (L_{t-1} - m_{t-1})\hat{f}_{t-1}(y_{t-1}^i) \right),
\end{multline}

\vspace{-4mm}

\noindent where $\hat{f}$ is the shifted objective as defined in \cref{eq:f_hat}.
\end{theorem}
\begin{proof}
We use the same notation as in the proof of \cref{thm:tracking:sector}. For some time index $t$, left and right multiply \cref{eq:thm2:lmi} by $\mathrm{vec}(\eta_t, u_t, \Delta \xi_t^\star, w_t)$ 
to obtain 
\begin{equation}
    \label{eq:proof_pain_thm_a}
    -\rho^2 \| \eta_t \|^2_{P_t} + \| \eta_{t+1} \|^2_{P_{t+1}} - \gamma_{\xi} \| \Delta \xi_t^\star \|^2 - \gamma_g \sum_{i=1}^p \| \Delta g_t(y_t^i) \|^2 + \sum_{i=1}^p \lambda_{p}^i (z_t^i)^\top \hat{M} z_t^i \leq 0.
\end{equation}
Next, multiply (\ref{eq:proof_pain_thm_a}) by $\rho^{2(k-t-1)}$ and sum up from $t=0$ to $t=k-1$ to obtain
\begin{multline}\label{eq:thm2:proof:2}
    \norm{\eta_{k}}^2_{P_k} \leq \rho^{2k} \norm{\eta_0}^2_{P_0} +  \sum_{t=0}^{k-1} \rho^{2(k-t-1)} \Biggl( \gamma_{\xi} \norm{\Delta \xi^\star_t}^2 + \gamma_g \sum_{i=1}^p \| \Delta g_t(y_t^i) \|^2 \\[-1ex]
    -  \sum_{i=1}^p \lambda_{p}^i (z_t^i)^\top \hat{M} z_t^i \Biggr).
\end{multline}
We now apply several bounds: first, apply $\underline{\lambda} \| \eta_i \|^2 \leq \norm{\eta_{i}}^2_{P_i} \leq \overline{\lambda} \| \eta_i \|^2$ $\forall i$ to \cref{eq:thm2:proof:2} and divide by $\underline{\lambda}$, yielding the constants $c_1$ and $c_2$ on the r.h.s. Second, use the fact that $\| \tilde{\xi}_k \| \leq \| \eta_k \|$. Third, assume w.l.o.g. that the IQC filter is initialized as $\zeta_0 = 0$, such that $\| \eta_0\|  = \| \mathrm{vec}(\tilde{\xi}_0, \zeta_0) \| = \| \tilde \xi_0 \|$. Thus, $\eta$ in \cref{eq:thm2:proof:2} can be replaced with $\tilde{\xi}$. Then, by a shift of the sum indices, we obtain
\begin{multline}
    \norm{\tilde{\xi}_k}^2 \leq c_1 \rho^{2k} \norm{\tilde{\xi}_0}^2 +  c_2 \sum_{t=1}^{k} \rho^{2(k-t)} \left( \gamma_{\xi} \norm{\Delta \xi^\star_{t-1}}^2 + \gamma_g \sum_{i=1}^p \| \Delta g_{t-1}(y_{t-1}^i) \|^2 \right) \\
    - c_2 \rho^{2(k-1)} \biggl( \sum_{i=1}^p \lambda_{p}^i \hspace{-5mm} \underbrace{\sum_{t=0}^{k-1} \rho^{-2t} (z_t^i)^\top \hat{M} z_t^i}_{\parbox{4cm}{\scriptsize \centering
        $\geq$ r.h.s.\ of \cref{eq:prop:off_by_1_iqc:ineq}\\
       $\quad$ evaluated at $y^i_t$
    }} \hspace{-5mm} \biggr)
\end{multline}
After plugging in the r.h.s. of \cref{eq:prop:off_by_1_iqc:ineq}, we pull $\rho^{2(k-1)}$ back into the sum, which yields \cref{eq:thm2:bound}.
\end{proof}

\cref{thm:tracking:variational} establishes a tracking bound that is more insightful than \cref{thm:tracking:sector} in the sense that it captures the dependence on multiple variation measures of the time-varying problem. 
We will show in the numerical examples that the use of variational IQCs will be particularly beneficial to obtain smaller worst-case decay rates $\rho$. 
Notably, the use of the variational IQC introduces an additional term that reflects changes in the function landscape, albeit in a form that is not directly interpretable. 
Nonetheless, \cref{thm:tracking:variational} provides valuable structural insights and forms the starting point for several corollaries.

\begin{corollary}\label{corr:constantSector}
Assume \cref{thm:tracking:variational} holds and that $m$ and $L$ are parameter-independent. W.l.o.g. define $y_{-1}^i=x_0^\star,\, i\in\mathbb{I}_p$ and $\hat{f}_{-1} \equiv 0$. Define $\gamma_f^i := \lambda_p^i \rho^2 (L - m)$. Then it holds
\begin{multline}\label{eq:corr1:bound}
    \norm{\xi_k - \xi_k^\star}^2 \leq c_1 \rho^{2k} \norm{\xi_0 - \xi_0^\star}^2 \\ + c_2 \sum_{t=0}^{k-1} \rho^{2(k-t-1)} \biggl( \gamma_{\xi} \norm{\Delta \xi^\star_{t}}^2 + \sum_{i=1}^p \gamma_g \| \Delta g_{t}(y_{t}^i) \|^2 - \gamma_f^i \Delta \hat{f}_{t-1}(y_{t-1}^i) \biggr)
\end{multline}
In particular, if $p=1$ we have
\begin{samepage}
\begin{multline}\label{eq:corr2:bound}
    \norm{\xi_k - \xi_k^\star}^2 \leq c_1 \rho^{2k} \norm{\xi_0 - \xi_0^\star}^2 + c_2 \sum_{t=0}^{k-1} \rho^{2(k-t-1)} \biggl( \gamma_{\xi} \norm{\xi^\star_{t+1} - \xi^\star_{t} }^2 \\ 
    + \gamma_g \| \nabla f_{t+1}(y_t) - \nabla f_t(y_t) \|^2 + \gamma_f \left( \hat{f}_{t}(y_{t-1}) - \hat{f}_{t-1}(y_{t-1}) \right) \biggr).
\end{multline}
\end{samepage}
\end{corollary}

In practice, one might have some knowledge about of the variational measures. For example, a common assumption is to have Lipschitz continuity of them w.r.t. $\theta$ \cite{simonetto_timeStructuredReview,simonetto_averagedOperators}, with which the following asymptotical bound can be established.

\begin{corollary}\label{corr:asymptotical}
    
    Let the assumptions of \cref{corr:constantSector} hold, and let \cref{assum:fixed_point_assum} hold with constant $U$. Suppose that the following variational measures are Lipschitz continuous in the parameter $\theta$
    \begin{samepage}
    \begin{align*}
        \norm{\nabla_x f(x, \theta_{k+1}) - \nabla_x f(x, \theta_{k})} &\leq  L_g \| \theta_{k+1} - \theta_{k} \|\quad \forall x\in\R^d\\
        | f(x, \theta_{k+1}) -  f(x, \theta_{k}) | &\leq  L_f \| \theta_{k+1} - \theta_{k} \|\quad \forall x\in\R^d\\
        |f(\text{\small $x^\star(\theta_{k+1})$}, \theta_{k+1}) - f(\text{\small $x^\star(\theta_k)$}, \theta_{k})| &\leq  L_{f^\star} \| \theta_{k+1} - \theta_{k} \|
    \end{align*}
    \end{samepage}
    for some $L_{g}, L_{f}, L_{f^\star} > 0$ for any $\theta_k,\theta_{k+1}\in\Theta$. Then we have
    \begin{equation}\label{eq:bound_asymptotical}
     \limsup_{k\rightarrow \infty} \norm{ \xi_k - \xi^\star_k} \leq \sqrt{c_2 \frac{(\gamma_\xi L_\xi^2 + \bar\gamma_g L_g^2) \nu_{\max}^2 + \bar{\gamma}_f  (L_f + L_{f^\star})\nu_{\max}}{1-\rho^2}},
    \end{equation}
    with $\bar\gamma_g \!:=\! p \gamma_g$, $\bar{\gamma}_f \!:=\! \sum_{i=1}^{p} \gamma_f^i$, $\nu_{\max} \!=\! \| \mathrm{vec}(\max\{ |\underline{\nu}_1|, |\overline{\nu}_1| \}, \dots, \max\{ |\underline{\nu}_{n_\theta}|, |\overline{\nu}_{n_\theta}| \}) \|$, and $L_\xi = \frac{L_g}{\underline{m}}\|U\|$.
    
\end{corollary}
\begin{proof}
    
    Lipschitz continuity of $\nabla_x f$ in $\theta$ together with uniform strong convexity of $f$ in $x$ implies $\| x^\star(\theta_{k+1}) - x^\star(\theta_k) \| \leq \frac{1}{\underline{m}} \| \nabla_x f(x^\star(\theta_{k}),\theta_{k+1}) - \nabla_x f(x^\star(\theta_{k}),\theta_k) \| \leq  \frac{L_g}{\underline{m}} \| \theta_{k+1} - \theta_{k} \|$. Using $\xi_k^\star = U x_k^\star$, we bound $\|\xi_{k+1}^\star - \xi_k^\star \| \leq L_\xi \| \theta_{k+1} - \theta_k\|$ with $L_\xi = \frac{L_g}{\underline{m}}\|U\|$. Now use the rate bound \cref{eq:param_rate_bound} for $\| \theta_{k+1} - \theta_k \| \leq \nu_{\max}$ and upper bound $\|\xi_{t+1}^\star - \xi_t^\star \|^2$, $\norm{\nabla f_{t+1}(y_t) - \nabla f_t(y_t)}^2$ and $| \hat f_{t}(y_{t-1}) - \hat f_{t-1}(y_{t-1})|$ in \cref{eq:corr2:bound} by $L_\xi^2 \nu_{\max}^2$, $L_g^2 \nu_{\max}^2$ and $(L_f + L_{f^\star}) \nu_{\max}$, respectively. Since $\rho \in (0,1)$, we may apply the geometric series formula to obtain \cref{eq:bound_asymptotical}.
\end{proof}

\begin{rem}\label{remark:bound_tradeoff}\label{rk:bound_tradeoff}
    Note that $\lambda_p$, $\gamma_\xi$ and $\gamma_g$ enter (\ref{eq:thm2:lmi}) affinely, and thus, they can be included as decision variables in an SDP. 
    The same holds for the condition number of $P(\theta)$, which enters \cref{eq:thm2:bound} indirectly through $c_1$ and $c_2$. We propose to trade off their values by solving a minimization problem
    \begin{align}\label{eq:SDP}
        \begin{aligned}
        \min_{P, \lambda_p, \gamma_{\xi}, \gamma_g, \sigma}&\quad  \sigma + k_1 \gamma_\xi + k_2 \gamma_g + k_3 \sum_{i=1}^{p}\lambda_p^i \\
        \mathrm{s.t.} &\quad (\ref{eq:thm2:lmi}), \\ & \quad  \lambda_p^i, \gamma_\xi, \gamma_g \geq 0, \quad \sigma \geq 1 \\
        &\quad P(\theta) \preceq \sigma I, \quad \begin{bmatrix}
        P(\theta) & I \\ I & \sigma I
        \end{bmatrix} \succeq 0 \quad \forall \theta \in \Theta
        \end{aligned}
    \end{align}
    for some non-negative weights $k_1,k_2,k_3 \geq 0$. The last two constraints are equivalent to $\frac{1}{\sigma} I \preceq P(\theta) \preceq \sigma I$, which bounds the condition number of $P(\theta)$ by $\sigma^2$. One may perform a bisection or a line search on $\rho$ to find the smallest bound \cref{eq:thm2:bound}. An optimal choice of the weights may be problem-specific, i.e. dependent on the foresight of the minimizer, gradient and function variation. 
\end{rem}

\begin{rem}
Typically, the values of stepsizes and other algorithm parameters are chosen based on $m$ and $L$. In practice, these quantities may only be approximately known and one advantage of the proposed robust control framework is that the effect of such misspecifications on the algorithm performance can naturally be analyzed. Given estimates $\hat m$ and $\hat L$ used for parameter tuning and error bounds with respect to the true values, we can model this mismatch as a parametric uncertainty, and resort to standard procedures in robust control (e.g. \cite{pfifer,padmanabhan}). We leave a detailed investigation of this and other potential strategies to future work.
\end{rem}

For \cref{thm:tracking:variational}, an iterate bound analogous to \cref{corr:iterates} can be derived along the same lines. We omit the explicit statement, as it follows by the same arguments.

Together, these results form a foundation for the numerical analysis in the next section, where we demonstrate the practical relevance of our approach and compare different algorithmic choices under time variation.

\section{Case Studies}\label{sec:case_studies} 
We now explore the practical implications of the theoretical results derived in this paper. The aim is to demonstrate the influence of algorithmic structures, rate bounds \cref{eq:param_rate_bound}, and type of IQC on the obtainable convergence rates and tracking bounds.

\subsection{Tracking Certificates}

Consider first an arbitrary objective $f \in \Smlt$ with a strong convexity parameter $m=1$ that remains constant, and a smoothness parameter $L_k$ that varies over time within an interval $\left[ \frac{4}{5}L_{\mathrm{nom}}, L_{\mathrm{nom}} \right]$ for some $L_{\mathrm{nom}} > \frac{5}{4}$, with a rate bound $|L_{k+1} - L_k| \leq \bar{\nu}$. We can model this within \cref{eq:min_fx_varying} by letting $\theta_k \triangleq L_k$. Note that $\nu_{\max} = \frac{1}{5}L_\mathrm{nom}$. Using \cref{thm:tracking:sector} and \cref{thm:tracking:variational}, we can investigate how the performance of optimization algorithms suffers from such variations and especially how the rate of change of $L_k$ impacts the achievable convergence rates. All algorithms used are tuned (or scheduled) with the varying smoothness parameter, and we plot our certificates over the (worst-case) function class condition ratio $L_\mathrm{nom}/m ~=~L_\mathrm{nom}$.

\subsubsection*{Obtaining Convergence Rates}
Let us first investigate the smallest achievable convergence rate $\rho$ in dependence of the algorithm, the choice of IQC, and the rate bound $\bar{\nu}$. We run a bisection on $\rho$ by solving the respective LMIs of \cref{thm:tracking:sector} and \cref{thm:tracking:variational}. We use a linear Lyapunov matrix parametrization $P(\theta) = P_0 + P_1 \theta$, with $P_0$ and $P_1$ being among the decision variables of the SDP. We note that our results did not improve further upon choosing a higher polynomial order. The LMIs are solved with a gridding approach \cite{pfifer} using \texttt{cvxpy} \cite{cvxpy} with solver \texttt{MOSEK}\footnote{The open-source implementation for all numerical examples presented in this paper can be accessed at: \faGithub\ \href{https://github.com/col-tasas/2024-tvopt-algorithm-analysis}{\texttt{https://github.com/col-tasas/2024-tvopt-algorithm-analysis}}.}.

\begin{figure}[ht]
	\centering
	\begin{subfigure}[t]{0.47\textwidth}
		\centering

\definecolor{mycolor1}{rgb}{0.8, 0.2, 0.247}   % red
\definecolor{mycolor2}{rgb}{0.267, 0.047, 0.329}   % dark purple
\definecolor{mycolor3}{rgb}{0.208, 0.373, 0.553}   % blue
\definecolor{mycolor4}{rgb}{0.133, 0.659, 0.518}   % teal
\definecolor{mycolor5}{rgb}{0.478, 0.824, 0.318}   % green
\definecolor{mycolor6}{rgb}{1.0, 0.6, 0.2}         % orange

\begin{tikzpicture}

\begin{axis}[%
width=1.8in,
height=1in,
at={(0.758in,0.481in)},
scale only axis,
xmode=log,
xmin=1,
grid style={dashed},
xmax=100,
xminorticks=true,
xlabel style={font=\color{white!15!black}},
xlabel={$L_\mathrm{nom}$},
ymin=0.00,
ymax=1.01,
ylabel style={font=\color{white!15!black}},
ylabel={$\rho$},
axis background/.style={fill=white},
xmajorgrids,
xminorgrids,
ymajorgrids,
legend style={at={(0.97,0.03)}, reverse legend, nodes={scale=0.75, transform shape}, anchor=south east, legend cell align=left, align=left, draw=white!15!black}]

\addplot [color=mycolor6, line width=1pt]
  table[row sep=crcr]{%
  1.251	0.112060546875\\
  2.33925569	0.4013671875\\
  4.37419438	0.62841796875\\
  8.17934379	0.7822265625\\
  15.29462549	0.87744140625\\
  28.59955208	0.93310546875\\
  53.47854904	0.9638671875\\
  100.0	0.980712890625\\
};
\addlegendentry{$L=L_{\mathrm{nom}}$ (static)}

\addplot [color=mycolor1, dashed, mark=+, mark options={solid, mycolor1}, line width=1pt]
  table[row sep=crcr]{%
  1.251	0.112060546875\\
  2.33925569	0.4013671875\\
  4.37419438	0.62841796875\\
  8.17934379	0.78515625\\
  15.29462549	0.884765625\\
  28.59955208	0.942626953125\\
  53.47854904	0.973388671875\\
  100.0	0.990234375\\
};
\addlegendentry{$\bar{\nu} = 0.05\bar{\nu}_{\max}$}

\addplot [color=mycolor3, dashed, mark=o, mark options={solid, mycolor3}, line width=1pt]
  table[row sep=crcr]{%
  1.251	0.11865234375\\
  2.33925569	0.431396484375\\
  4.37419438	0.68408203125\\
  8.17934379	0.859130859375\\
  15.29462549	0.971923828125\\
};
\addlegendentry{$\bar{\nu} = 0.5\bar{\nu}_{\max}$}

\addplot [color=mycolor4, dashed, mark=x, mark options={solid, mycolor4}, line width=1pt]
  table[row sep=crcr]{%
  1.251	0.1259765625\\
  2.33925569	0.466552734375\\
  4.37419438	0.75\\
  8.17934379	0.9521484375\\
};
\addlegendentry{$\bar{\nu} = \bar{\nu}_{\max}$}

\end{axis}
\end{tikzpicture}%
		\caption{Influence of higher parameter rate bounds $\bar{\nu}$ (one-step GD).}
		\label{fig:pointwise:nu}
	\end{subfigure}
	\hfill
	\begin{subfigure}[t]{0.47\textwidth}
		\centering

\definecolor{mycolor1}{rgb}{0.8, 0.2, 0.247}   % red
\definecolor{mycolor2}{rgb}{0.267, 0.047, 0.329}   % dark purple
\definecolor{mycolor3}{rgb}{0.208, 0.373, 0.553}   % blue
\definecolor{mycolor4}{rgb}{0.133, 0.659, 0.518}   % teal
\definecolor{mycolor5}{rgb}{0.478, 0.824, 0.318}   % green
\definecolor{mycolor6}{rgb}{1.0, 0.6, 0.2}         % orange

\begin{tikzpicture}

\begin{axis}[%
width=1.8in,
height=1in,
at={(0.758in,0.481in)},
scale only axis,
xmode=log,
xmin=1,
xmax=100,
grid style={dashed},
xminorticks=true,
xlabel style={font=\color{white!15!black}},
xlabel={$L_\mathrm{nom}$},
ymin=0.00,
ymax=1.01,
ylabel style={font=\color{white!15!black}},
ylabel={$\rho$},
axis background/.style={fill=white},
xmajorgrids,
xminorgrids,
ymajorgrids,
legend style={at={(0.97,0.03)}, reverse legend, nodes={scale=0.65, transform shape}, anchor=south east, legend cell align=left, align=left, draw=white!15!black}]

\addplot [color=mycolor1, dashed, mark=+, mark options={solid, mycolor1}, line width=1pt]
  table[row sep=crcr]{%
  1.251	0.112060546875\\
  2.33925569	0.4013671875\\
  4.37419438	0.62841796875\\
  8.17934379	0.78515625\\
  15.29462549	0.884765625\\
  28.59955208	0.942626953125\\
  53.47854904	0.973388671875\\
  100.0	0.990234375\\
};
\addlegendentry{GD}

\addplot [color=mycolor2, dashed, mark=triangle, mark options={solid, mycolor2}, line width=1pt]
  table[row sep=crcr]{%
1.251	        0.01318359375\\
2.03550964	    0.116455078125\\
3.31199001	    0.287841796875\\
5.38895892	    0.472412109375\\
8.76840756	    0.6328125\\
14.2671288	    0.755859375\\
23.2141313	    0.843017578125\\
37.77185302	    0.904541015625\\
61.4588098	    0.944091796875\\
100.0	        0.9697265625\\
};
\addlegendentry{2-Step GD}

\addplot [color=mycolor4, dashed, mark=square, mark options={solid, mycolor4}, line width=1pt]
  table[row sep=crcr]{%
1.251	        0.000732421875\\
2.03550964	    0.005126953125\\
3.31199001	    0.044677734375\\
5.38895892	    0.153076171875\\
8.76840756	    0.318603515625\\
14.2671288	    0.495849609375\\
23.2141313	    0.651123046875\\
37.77185302	    0.77197265625\\
61.4588098	    0.85693359375\\
100.0	        0.913330078125\\
};
\addlegendentry{5-Step GD}

% \addplot [color=mycolor4, dashed, mark=o, mark options={solid, mycolor4}, line width=1pt]
%   table[row sep=crcr]{%
%   1.25	0.142822265625\\
%   2.03406326242116	0.37353515625\\
%   3.30993068442511	0.533335546875\\
%   5.38608672507971	0.6489\\
%   8.76451290855913	0.7338\\
%   14.2620589762529	0.8000\\
%   23.2079441680639	0.8562\\
%   37.7651413028641	0.9074\\
%   61.4533492194397	0.9484\\
%   100	0.9814\\
% };
% \addlegendentry{Nesterov}

% \addplot [color=mycolor5, dashed, mark=x, mark options={solid, mycolor5}, line width=1pt]
%   table[row sep=crcr]{%
% 1.251	        0.106201171875\\
% 2.03550964	    0.299560546875\\
% 3.31199001	    0.451171875\\
% 5.38895892	    0.5712890625\\
% 8.76840756	    0.665771484375\\
% 14.2671288	    0.744873046875\\
% 23.2141313	    0.812255859375\\
% 37.77185302	    0.8671875\\
% 61.4588098	    0.91259765625\\
% 100.0	        0.94921875\\
% };
% \addlegendentry{Triple Mom.}
4\end{axis}
\end{tikzpicture}%
		\caption{Influence of higher number of algorithm evaluations ($\bar{\nu}=0.05\nu_{\max}$).}
		\label{fig:pointwise:algos}
	\end{subfigure}
    \vspace{-1em}
	\caption{Influence of rate bounds $\overline{\nu}$ and number of algorithm evaluations, for GD.}
	\label{fig:pointwise}
	\vspace{-3mm}
\end{figure}
Fig.~\ref{fig:pointwise} shows the convergence rates obtained by gradient descent (GD) ($\bar{A} = I_d,\, \bar{B}(\theta_k) = - \tfrac{1}{L_k} I_d,\, \bar{C} = I_d$) for different settings. First, we fix the algorithm and let $\overline{\nu}$ increase (Fig.~\ref{fig:pointwise:nu}), and, second, we fix $\overline{\nu}$ and let $p$ increase (Fig.~\ref{fig:pointwise:algos}), i.e., we consider multi-step GD.
A few observations are in order. First, the convergence rates $\rho$ become worse for high condition ratios $L_\mathrm{nom}$, which is a well known phenomenon for strongly convex and smooth objectives. Second, observe from Fig.~\ref{fig:pointwise:nu} that higher rate bounds $\overline{\nu}$ have a detrimental effect on $\rho$, while a vanishing rate bound recovers the same $\rho$ that is known from static optimization. Third, we can observe the beneficial effect of increasing the number of steps $p$ in Fig.~\ref{fig:pointwise:algos}.

We note that for GD, as well as for its $p$-step versions, the minimum value of $\rho$ was the same for \cref{thm:tracking:sector} and \cref{thm:tracking:variational}. However, we show in Fig.~\ref{fig:iqcComparison} that the use of a variational IQC is indeed beneficial compared to a pointwise IQC for analyzing accelerated first-order methods. Keeping $\overline{\nu} = 0.05 \nu_{\max}$ fixed, we compare GD with Nesterov's method (\cref{exmp:accelerated_algo} with $\alpha_k = \frac{1}{L_k}$, $\beta_k = \frac{\sqrt{L_k}-\sqrt{m_k}}{\sqrt{L_k}+\sqrt{m_k}}$, $\gamma_k=\beta_k$) and Triple Momentum ($\alpha_k = \frac{2\sqrt{L_k}-\sqrt{m_k}}{L_k\sqrt{L_k}}$, $\beta_k = \frac{(\sqrt{L_k}-\sqrt{m_k})^2}{\sqrt{L_k}(\sqrt{L_k}+\sqrt{m_k})}$, $\gamma_k=\frac{(\sqrt{L_k}-\sqrt{m_k})^2}{(2\sqrt{L_k}-\sqrt{m_k})(\sqrt{L_k}+\sqrt{m_k})}$), where we observe clear improvements in Fig.~\ref{fig:iqcComparison:variational} over Fig.~\ref{fig:iqcComparison:pointwise}, thus highlighting the advantage of the reduced conservatism of dynamic IQCs. 

\begin{figure}
	\centering
	\begin{subfigure}[t]{0.47\textwidth}
		\centering

\definecolor{mycolor1}{rgb}{0.8, 0.2, 0.247}   % red
\definecolor{mycolor2}{rgb}{0.267, 0.047, 0.329}   % dark purple
\definecolor{mycolor3}{rgb}{0.208, 0.373, 0.553}   % blue
\definecolor{mycolor4}{rgb}{0.133, 0.659, 0.518}   % teal
\definecolor{mycolor5}{rgb}{0.478, 0.824, 0.318}   % green
\definecolor{mycolor6}{rgb}{1.0, 0.6, 0.2}         % orange

\begin{tikzpicture}

\begin{axis}[%
width=1.8in,
height=1in,
at={(0.758in,0.481in)},
scale only axis,
xmode=log,
xmin=1,
grid style={dashed},
xmax=100,
xminorticks=true,
xlabel style={font=\color{white!15!black}},
xlabel={$L_\mathrm{nom}$},
ymin=0.11,
ymax=1.01,
ylabel style={font=\color{white!15!black}},
ylabel={$\rho$},
axis background/.style={fill=white},
xmajorgrids,
xminorgrids,
ymajorgrids,
legend style={at={(0.97,0.03)}, reverse legend, nodes={scale=0.75, transform shape}, anchor=south east, legend cell align=left, align=left, draw=white!15!black}]

\addplot [color=mycolor1, dashed, mark=+, mark options={solid, mycolor1}, line width=1pt]
  table[row sep=crcr]{%
  1.251	0.112060546875\\
  2.33925569	0.4013671875\\
  4.37419438	0.62841796875\\
  8.17934379	0.78515625\\
  15.29462549	0.884765625\\
  28.59955208	0.942626953125\\
  53.47854904	0.973388671875\\
  100.0	0.990234375\\
};
\addlegendentry{GD}

\addplot [color=mycolor3, dashed, mark=o, mark options={solid, mycolor3}, line width=1pt]
  table[row sep=crcr]{%
  1.251	0.150146484375\\
  2.03550964314036	0.42333984375\\
  3.31199001384285	0.63720703125\\
  5.38895892179193	0.8056640625\\
  8.76840755539149	0.936767578125\\
};
\addlegendentry{NM}

\addplot [color=mycolor4, dashed, mark=x, mark options={solid, mycolor4}, line width=1pt]
  table[row sep=crcr]{%
1.251	        0.1171875\\
2.03550964	    0.37939453125\\
3.31199001	    0.62548828125\\
5.38895892	    0.841552734375\\
};
\addlegendentry{TM}

\end{axis}
\end{tikzpicture}%
		\caption{Using \cref{thm:tracking:sector} ($\bar{\nu}=0.05\nu_{\max}$).}
		\label{fig:iqcComparison:pointwise}
	\end{subfigure}
	\hfill
	\begin{subfigure}[t]{0.47\textwidth}
		\centering

\definecolor{mycolor1}{rgb}{0.8, 0.2, 0.247}   % red
\definecolor{mycolor2}{rgb}{0.267, 0.047, 0.329}   % dark purple
\definecolor{mycolor3}{rgb}{0.208, 0.373, 0.553}   % blue
\definecolor{mycolor4}{rgb}{0.133, 0.659, 0.518}   % teal
\definecolor{mycolor5}{rgb}{0.478, 0.824, 0.318}   % green
\definecolor{mycolor6}{rgb}{1.0, 0.6, 0.2}         % orange

\begin{tikzpicture}

\begin{axis}[%
width=1.8in,
height=1in,
at={(0.758in,0.481in)},
scale only axis,
xmode=log,
grid style={dashed},
xmin=1,
xmax=100,
xminorticks=true,
xlabel style={font=\color{white!15!black}},
xlabel={$L_\mathrm{nom}$},
ymin=0.11,
ymax=1.01,
ylabel style={font=\color{white!15!black}},
ylabel={$\rho$},
axis background/.style={fill=white},
xmajorgrids,
xminorgrids,
ymajorgrids,
legend style={at={(0.97,0.03)}, reverse legend, nodes={scale=0.75, transform shape}, anchor=south east, legend cell align=left, align=left, draw=white!15!black}]

\addplot [color=mycolor1, dashed, mark=+, mark options={solid, mycolor1}, line width=1pt]
  table[row sep=crcr]{%
 1.251	0.112060546875\\
  2.33925569	0.4013671875\\
  4.37419438	0.62841796875\\
  8.17934379	0.78515625\\
  15.29462549	0.884765625\\
  28.59955208	0.942626953125\\
  53.47854904	0.973388671875\\
  100.0	0.990234375\\
};
\addlegendentry{GD}

\addplot [color=mycolor3, dashed, mark=o, mark options={solid, mycolor3}, line width=1pt]
  table[row sep=crcr]{%
1.25	0.142822265625\\
2.03406326242116	0.37353515625\\
3.30993068442511	0.533335546875\\
5.38608672507971	0.6489\\
8.76451290855913	0.7338\\
14.2620589762529	0.8000\\
23.2079441680639	0.8562\\
37.7651413028641	0.9074\\
61.4533492194397	0.9484\\
100	0.9814\\
};
\addlegendentry{NM}

\addplot [color=mycolor4, dashed, mark=x, mark options={solid, mycolor4}, line width=1pt]
  table[row sep=crcr]{%
1.251	        0.106201171875\\
2.03550964	    0.299560546875\\
3.31199001	    0.451171875\\
5.38895892	    0.5712890625\\
8.76840756	    0.665771484375\\
14.2671288	    0.744873046875\\
23.2141313	    0.812255859375\\
37.77185302	    0.8671875\\
61.4588098	    0.91259765625\\
100.0	        0.94921875\\
};
\addlegendentry{TM}

\end{axis}
\end{tikzpicture}%
		\caption{Using \cref{thm:tracking:variational} ($\bar{\nu}=0.05\nu_{\max}$).}
		\label{fig:iqcComparison:variational}
	\end{subfigure}
    \vspace{-1em}
	\caption{Comparison of convergence rates obtained by \cref{thm:tracking:sector} and \cref{thm:tracking:variational}.}
	\label{fig:iqcComparison}
\end{figure}

\begin{figure}
  \centering
  \begin{subfigure}[t]{0.325\textwidth}
    \centering
    \definecolor{mycolor1}{rgb}{0.8, 0.2, 0.247}   % red
\definecolor{mycolor3}{rgb}{0.208, 0.373, 0.553}   % blue
\definecolor{mycolor4}{rgb}{0.133, 0.659, 0.518}   % teal

\begin{tikzpicture}

\begin{axis}[%
width=1in,
height=1in,
at={(0.758in,0.481in)},
scale only axis,
xmode=log,
xmin=1,
xmax=100,
xminorticks=true,
xlabel style={font=\color{white!15!black}},
xlabel={$L_\mathrm{nom}$},
ymin=0.00,
ymax=2.6,
ylabel style={font=\color{white!15!black}, at={(axis description cs:-0.15,0.5)},
  anchor=south, yshift=-0.5em},
ylabel={$\lambda_p$},
axis background/.style={fill=white},
xmajorgrids,
% xminorgrids,
ymajorgrids,
grid style={dashed},
% xtick={1,10},
% tick label style={/pgf/number format/fixed}, % optional
legend style={at={(0.97,0.97)}, reverse legend, nodes={scale=0.65, transform shape}, anchor=north east, legend cell align=left, align=left, draw=white!15!black}]

\addplot [color=mycolor1, dashed, mark=+, mark options={solid, mycolor1}, line width=1pt]
  table[row sep=crcr]{%
1.251       0.0023609 \\
2.03550964    1.19591119\\
3.31199001    0.31923108\\
5.38895892    0.16870054\\
8.76840756    0.12326183\\
14.2671288    0.0825301 \\
23.2141313    0.06672526\\
37.77185302   0.04982391\\
61.4588098    0.03134629\\
100.          0.01953723\\
};
\addlegendentry{GD}

\addplot [color=mycolor3, dashed, mark=o, mark options={solid, mycolor3}, line width=1pt]
  table[row sep=crcr]{%
1.251     0.004255  \\
2.03550964 1.8359483 \\
3.31199001 0.55268209\\
5.38895892 0.22741888\\
8.76840756  0.11812157\\
14.2671288  0.06937018\\
23.2141313  0.03922371\\
37.77185302 0.02168277\\
61.4588098  0.01292884\\
100.         0.00730093\\
};
\addlegendentry{NM}

\addplot [color=mycolor4, dashed, mark=x, mark options={solid, mycolor4}, line width=1pt]
  table[row sep=crcr]{%
1.251        0.0042558 \\
2.03550964   2.5466957 \\
3.31199001   1.08528372\\
5.38895892   0.591041187\\
8.76840756   0.352766739\\
14.2671288   0.219299538\\
23.2141313   0.130538631\\
37.77185302  0.076729922\\
61.4588098   0.044630157\\
100.         0.024658646\\
};
\addlegendentry{TM}

\end{axis}
\end{tikzpicture}
    \caption{Sensitivity $\lambda_p$ $(p=1)$.}
  \end{subfigure}
  \begin{subfigure}[t]{0.325\textwidth}
    \centering
    \definecolor{mycolor1}{rgb}{0.8, 0.2, 0.247}   % red
\definecolor{mycolor3}{rgb}{0.208, 0.373, 0.553}   % blue
\definecolor{mycolor4}{rgb}{0.133, 0.659, 0.518}   % teal

\begin{tikzpicture}

\begin{axis}[%
width=1in,
height=1in,
at={(0.758in,0.481in)},
scale only axis,
xmode=log,
xmin=1,
xmax=100,
xminorticks=true,
xlabel style={font=\color{white!15!black}},
xlabel={$L_\mathrm{nom}$},
ymin=0.00,
ymax=1000,
ylabel style={font=\color{white!15!black}, yshift=-1em},
ylabel={$\gamma_\xi$},
axis background/.style={fill=white},
xmajorgrids,
% xminorgrids,
ymajorgrids,
grid style={dashed},
legend style={at={(0.03,0.97)}, reverse legend, nodes={scale=0.65, transform shape}, anchor=north west, legend cell align=left, align=left, draw=white!15!black}
]

\addplot [color=mycolor1, dashed, mark=+, mark options={solid, mycolor1}, line width=1pt]
  table[row sep=crcr]{%
1.251        77.011751581\\
2.03550964   52.28408018\\
3.31199001   45.77647871\\
5.38895892   36.174116649\\
8.76840756   27.933435879\\
14.2671288   29.31092003\\
23.2141313   21.8428656\\
37.77185302  11.75704728\\
61.4588098   10.22457897\\
100.         5.0276628007\\
};
\addlegendentry{GD}

\addplot [color=mycolor3, dashed, mark=o, mark options={solid, mycolor3}, line width=1pt]
  table[row sep=crcr]{%
1.251        117.518815534\\
2.03550964    38.11954495\\
3.31199001    76.867013717\\
5.38895892   101.823848642\\
8.76840756   124.158050558\\
14.2671288   144.624527822\\
23.2141313   193.514343314\\
37.77185302  285.8710374753\\
61.4588098   393.132223286\\
100.         616.9727215245\\
};
\addlegendentry{NM}

\addplot [color=mycolor4, dashed, mark=x, mark options={solid, mycolor4}, line width=1pt]
  table[row sep=crcr]{%
1.251        117.518645\\
2.03550964   166.011176\\
3.31199001   237.971128\\
% 5.38895892   585.480606\\
8.76840756   270.642763\\
14.2671288   326.054272\\
23.2141313   594.51809272\\
37.77185302  876.884835\\
% 61.4588098   1510.298003\\
100.         899.32609788\\
};
\addlegendentry{TM}

\end{axis}
\end{tikzpicture}
    \caption{Sensitivity $\gamma_\xi$.}
  \end{subfigure}
  \begin{subfigure}[t]{0.325\textwidth}
    \centering
    \definecolor{mycolor1}{rgb}{0.8, 0.2, 0.247}   % red
\definecolor{mycolor3}{rgb}{0.208, 0.373, 0.553}   % blue
\definecolor{mycolor4}{rgb}{0.133, 0.659, 0.518}   % teal

\begin{tikzpicture}

\begin{axis}[%
width=1in,
height=1in,
at={(0.758in,0.481in)},
scale only axis,
xmode=log,
xmin=1,
xmax=100,
xminorticks=true,
xlabel style={font=\color{white!15!black}},
xlabel={$L_\mathrm{nom}$},
ymin=0.00,
ymax=200,
ylabel style={font=\color{white!15!black},
  anchor=south, yshift=-0.5em},
ylabel={$\gamma_g$},
axis background/.style={fill=white},
xmajorgrids,
% xminorgrids,
ymajorgrids,
grid style={dashed},
legend style={at={(0.97,0.97)}, reverse legend, nodes={scale=0.65, transform shape}, anchor=north east, legend cell align=left, align=left, draw=white!15!black}
]

\addplot [color=mycolor1, dashed, mark=+, mark options={solid, mycolor1}, line width=1pt]
  table[row sep=crcr]{%
1.251        0.720614457\\
2.03550964   51.79927294\\
3.31199001   45.9890178\\
5.38895892   36.26051913\\
8.76840756   28.354379\\
14.2671288   29.921558\\
23.2141313   22.8850304\\
37.77185302  13.13835192\\
61.4588098   11.8015029\\
100.         7.55258327\\
};
\addlegendentry{GD}

\addplot [color=mycolor3, dashed, mark=o, mark options={solid, mycolor3}, line width=1pt]
  table[row sep=crcr]{%
1.251        0.59415971\\
2.03550964   11.948019\\
3.31199001   7.5027062\\
5.38895892   3.4112819\\
8.76840756   1.5787828\\
14.2671288   0.96857663\\
23.2141313   0.85406479\\
37.77185302  0.85725681\\
61.4588098   0.87104368\\
100.         0.88406697\\
};
\addlegendentry{NM}

\addplot [color=mycolor4, dashed, mark=x, mark options={solid, mycolor4}, line width=1pt]
  table[row sep=crcr]{%
1.251        0.594159308\\
2.03550964   146.9681832\\
3.31199001   173.599474\\
% 5.38895892   327.462554\\
8.76840756   102.5464942\\
14.2671288   67.423569476\\
23.2141313   69.587101255\\
37.77185302  61.79968246\\
61.4588098   64.5198716\\
100.         25.613581407\\
};
\addlegendentry{TM}

\end{axis}
\end{tikzpicture}
    \caption{Sensitivity $\gamma_g$.}
  \end{subfigure}
  \caption{Sensitivity of algorithms to variational measures. ($\bar{\nu}=0.05\nu_{\max}$)}
  \label{fig:sensitivity_analysis}
  \vspace{-3mm}
\end{figure}

\subsubsection*{Trade-offs of Tracking Certificates}

We note that although Fig.~\ref{fig:iqcComparison:variational} may suggest that accelerated algorithms like NM and TM converge faster (in worst case), it is important to remember that the tracking performance depends also on the variational measures, which are amplified by the sensitivities $\lambda_p$, $\gamma_\xi$, and $\gamma_g$. Fig.~\ref{fig:sensitivity_analysis} complements Fig.~\ref{fig:iqcComparison:variational} by plotting the values of $\lambda_p$, $\gamma_g$, and $\gamma_\xi$ which we obtain by solving \cref{eq:SDP}, here, for example, for $k_1,k_2,k_3=1$. Observe that TM, which has the best convergence rates, is significantly more sensitive to all variational measures, while GD is rather insensitive in comparison. This seems to resonate with the convergence-robustness trade-off that was shown for accelerated gradient methods in static optimization \cite{jovanovic_robustness}, and may explain why accelerated methods can perform worse in time-varying settings \cite{dallAnese_informationStreams}.

Note that for each value of $\rho$, we obtain one set of sensitivities $(\lambda_p, \gamma_\xi, \gamma_g)$ and constants $c_1,c_2$, and therefore one tracking bound. The bound that we obtain by solving \cref{eq:SDP} at the minimal $\rho$ is not necessarily the tightest. In fact, we can visualize a trade-off of $\rho$ vs. $(\sigma, \lambda_p, \gamma_\xi, \gamma_g)$, see Fig.~\ref{fig:trade_off:rho}. Recall that $\sigma^2$ upper bounds the condition number of $P(\theta)$, which influences the constants $c_1, c_2$. We may control the magnitude of each sensitivity via their respective weighting $k_1,k_2,k_3$ in \cref{eq:SDP}, however, suppressing the sensitivities may in turn lead to bad condition numbers of $P(\theta)$. We visualize the $\sigma^2$ vs. $(\lambda_p, \gamma_\xi, \gamma_g)$ trade-off in Fig. \ref{fig:trade_off:sigma}. We note that the same trade-offs can also be observed among the sensitivities themselves when weighting individual $k_i$ higher than the other, which is omitted here. Ultimately, obtaining the parameters that yield the best tracking bound is a non-trivial task. One might attempt to obtain a Pareto-optimal set of parameters by carefully balancing the values of $\rho$ and $k_1,k_2,k_3$, e.g. when prior knowledge of $\Delta \xi_k^\star$, $\Delta g_k$, $\Delta f_k$ is available.

\begin{figure}
	\centering
	\begin{subfigure}[t]{0.47\textwidth}
		\centering
\definecolor{mycolor1}{rgb}{0.8, 0.2, 0.247}   % red
\definecolor{mycolor2}{rgb}{0.267, 0.047, 0.329}   % dark purple
\definecolor{mycolor3}{rgb}{0.208, 0.373, 0.553}   % blue
\definecolor{mycolor4}{rgb}{0.133, 0.659, 0.518}   % teal
\definecolor{mycolor5}{rgb}{0.478, 0.824, 0.318}   % green
\definecolor{mycolor6}{rgb}{1.0, 0.6, 0.2}        
\begin{tikzpicture}

\begin{axis}[
  width=1.6in,
  height=1.1in,
  at={(0.758in,0.523in)},
  scale only axis,
  xmin=0.3,
  xmax=1,
  xminorticks=true,
  xlabel style={font=\color{white!15!black}},
  xlabel={$\rho$},
  ymin=0.00,
  ymax=26,
  yminorticks=true,
  ylabel={\footnotesize $\sigma$, $\gamma_\xi, \gamma_g, \lambda_p$},      % <<< right y-axis label
  ylabel style={font=\color{black}, yshift=0em},
  xlabel={$\rho$},
  xlabel style={font=\color{black}, yshift=-0.45em},  
  yticklabel pos=left,
  xmajorgrids,
  ymajorgrids,
  yminorgrids,
  grid style={dashed},   
  xtick={0.2,0.4,0.6,0.8,1,0.34},
  xticklabels={$0.2$, $\,$, $0.6$,$0.8$,$1$,$\textcolor{gray}{\rho_{\min}}$},
  xticklabel style={font=\footnotesize},
  legend style={at={(0.97,0.97)}, anchor=north east, legend cell align=left, nodes={scale=0.7, transform shape}}
]

\addplot [color=mycolor1, dashdotted, line width=1.2pt, mark=triangle, mark options={solid, mycolor1}]
  table[row sep=crcr]{%
0.35	61.5617036177269\\
0.4	22.5632122402881\\
0.45	14.2003423933622\\
0.5	10.4781270610043\\
0.55	8.38499383671709\\
0.6	7.06463610024617\\
0.65	6.26327333052195\\
0.7	5.61907556791067\\
0.75	5.14988191272652\\
0.8	4.79732280689582\\
0.85	4.52419982096216\\
0.9	4.30281587778014\\
0.95	4.11121378064539\\
1	3.93693831195623\\
};
\addlegendentry{$\sigma$}

\addplot [color=mycolor2, dashed, line width=1.2pt, mark=diamond, mark options={solid, mycolor2}]
  table[row sep=crcr]{%
0.35	25.6975280320547\\
0.4	8.62690397110993\\
0.45	4.97616887427592\\
0.5	3.28708968527981\\
0.55	2.30680126754184\\
0.6	1.68104995717341\\
0.65	1.29904725466851\\
0.7	1.24472679607495\\
0.75	1.20523250548008\\
0.8	1.1699870931456\\
0.85	1.13786274138816\\
0.9	1.111827093615\\
0.95	1.09497356559948\\
1	1.08679110705913\\
};
\addlegendentry{$\gamma_\xi$}

\addplot [color=mycolor3, dotted, line width=1.2pt, mark=square, mark options={solid, mycolor3}]
  table[row sep=crcr]{%
0.35	25.6334198180948\\
0.4	8.95375342890766\\
0.45	5.61012301714357\\
0.5	4.12254935241141\\
0.55	3.2724802033606\\
0.6	2.70772004161351\\
0.65	2.21374651904923\\
0.7	1.78586327911218\\
0.75	1.47417510257955\\
0.8	1.24757709171942\\
0.85	1.08450765382334\\
0.9	0.970843393220161\\
0.95	0.896638968200513\\
1	0.851289895830653\\
};
\addlegendentry{$\gamma_g$}

\addplot [color=mycolor4, dashdotted, line width=1.2pt, mark=triangle, mark options={solid, mycolor4}]
  table[row sep=crcr]{%
0.35	0.944101686878973\\
0.4	1.29092537665433\\
0.45	1.39311653295509\\
0.5	1.43965337424826\\
0.55	1.46395396666796\\
0.6	1.47600816282426\\
0.65	1.44919099072124\\
0.7	1.3571700009345\\
0.75	1.25962933596126\\
0.8	1.16362843514627\\
0.85	1.07115531744657\\
0.9	0.9819207238356\\
0.95	0.896174759940087\\
1	0.816005615951696\\
};
\addlegendentry{$\lambda_p$}

\addplot [gray, thick] coordinates {(0.335,0) (0.34,26)};

\end{axis}

\end{tikzpicture}%
		\caption{Trade-off between convergence rate and sensitivities/condition number, for fixed weightings $k_1,k_2,k_3=1$.}
		\label{fig:trade_off:rho}
	\end{subfigure}
	\hfill
	\begin{subfigure}[t]{0.47\textwidth}
		\centering
\definecolor{mycolor1}{rgb}{0.8, 0.2, 0.247}   % red
\definecolor{mycolor2}{rgb}{0.267, 0.047, 0.329}   % dark purple
\definecolor{mycolor3}{rgb}{0.208, 0.373, 0.553}   % blue
\definecolor{mycolor4}{rgb}{0.133, 0.659, 0.518}   % teal
\definecolor{mycolor5}{rgb}{0.478, 0.824, 0.318}   % green
\definecolor{mycolor6}{rgb}{1.0, 0.6, 0.2}        
\begin{tikzpicture}

\begin{axis}[
  width=1.6in,
  height=1.1in,
  at={(0.758in,0.523in)},
  scale only axis,
  xmin=0,
  xmax=110,
  separate axis lines,
  ymode=log,
  ymin=0.05,
  ymax=150,
  yminorticks=true,
  ylabel={\footnotesize $\gamma_\xi, \gamma_g, \lambda_p$},      % <<< right y-axis label
  ylabel style={font=\color{black}, yshift=0.9em},
  axis y line*=right,          % <<< right axis line
  axis x line*=none,           % <<< DON'T redraw x-axis
  yticklabel pos=right,
  xmajorgrids,
  ymajorgrids,
  % yminorgrids,
  grid style={dashed},
  ytick={0.1,1,10,100},       
  axis x line*=top,
  xtick=\empty,                            % ticks at 10^0 and 10^2
  yticklabels={{$0.1$},{$1$},{$10$},{$10^2$}},
  legend style={at={(0.97,0.47)}, fill=white, anchor=north east, legend cell align=left, nodes={scale=0.6, transform shape}}
]
\addlegendimage{empty legend}
\addlegendentry{}
\addlegendimage{empty legend}
\addlegendentry{}

\addplot [color=mycolor2, dashed, line width=1.2pt, mark=diamond, mark options={solid, mycolor2}]
  table[row sep=crcr]{%
1	0.556446857469742\\
23.1111111111111	0.209466535499259\\
45.2222222222222	0.145219773995414\\
67.3333333333333	0.124451721259083\\
89.4444444444444	0.0936795202955725\\
111.555555555556	0.0875966810148187\\
133.666666666667	0.0830648009318011\\
155.777777777778	0.0829046239696487\\
177.888888888889	0.078553978001286\\
200	0.0766348324900386\\
};
\addlegendentry{$\gamma_\xi$}

\addplot [color=mycolor3, dotted, line width=1.2pt, mark=square, mark options={solid, mycolor3}]
  table[row sep=crcr]{%
1	69.3990609890014\\
23.1111111111111	5.83051409588952\\
45.2222222222222	4.44927996326799\\
67.3333333333333	3.31380666039053\\
89.4444444444444	3.83344037500161\\
111.555555555556	3.16061191078473\\
133.666666666667	2.68006857488648\\
155.777777777778	2.11477806368691\\
177.888888888889	1.92326749167847\\
200	1.67414618205868\\
};
\addlegendentry{$\gamma_g$}

\addplot [color=mycolor4, dashdotted, line width=1.2pt, mark=triangle, mark options={solid, mycolor4}]
  table[row sep=crcr]{%
1	69.4371473085098\\
23.1111111111111	5.92262873901018\\
45.2222222222222	4.51236295093574\\
67.3333333333333	3.36601261719795\\
89.4444444444444	3.87528475309391\\
111.555555555556	3.19911424640859\\
133.666666666667	2.71607697254613\\
155.777777777778	2.15012806935029\\
177.888888888889	1.95654235697202\\
200	1.70579863919184\\
};
\addlegendentry{$\lambda_p$}

\end{axis}

\begin{axis}[
  width=1.6in,
  height=1.1in,
  at={(0.758in,0.523in)},
  scale only axis,
  xmin=0,
  xmax=110,
  xlabel={$k_1, k_2, k_3$},
  separate axis lines,
  ymode=log,
  ymin=1e4,
  ymax=1e6,
  yminorticks=true,
  ylabel={\small $\sigma^2$ (cond. of $P$)}, % <<< left y-axis label
  ylabel style={, yshift=-0.7em},
  axis y line*=left,          % <<< left axis line
  axis x line*=bottom,        % <<< shared x-axis
  xmajorgrids,
  ymajorgrids,
  % yminorgrids,
  grid style={dashed},
  legend style={at={(0.97,0.48)}, fill=none, anchor=north east, legend cell align=left,draw=none,  nodes={scale=0.6, transform shape}}
]
\addplot [color=mycolor1, line width=1.2pt, mark=o, mark options={solid, mycolor1}]
  table[row sep=crcr]{%
1                       22882.2784845198\\
23.1111111111111        76449.97809592652\\
45.2222222222222        169686.8720280791\\
67.3333333333333        209912.06698059\\
89.4444444444444        487021.8151602564\\
111.555555555556        517263.342795655\\
133.666666666667        536393.0786653651\\
155.777777777778        458725.3932567435\\
177.888888888889        495837.39608618716\\
200                     477924.83954424417\\
};
\addlegendentry{$\sigma^2$}
\end{axis}

\end{tikzpicture}%
		\caption{Trade-off between condition of $P(\theta)$ and sensitivities, for fixed $\rho=\rho_\mathrm{min}$.}
		\label{fig:trade_off:sigma}
	\end{subfigure}
    \vspace{-1em}
	\caption{Trade-off in obtaining tracking certificates from \cref{eq:SDP}. Numeric study for GD, $L_{\mathrm{nom}} = 2$, $\bar{\nu} = 0.05\nu_{\max}$.}
	\label{fig:trade_off}
	\vspace{-3mm}
\end{figure}

\subsection{Tracking Bounds}

We now visualize our bounds on an instance of a time-varying optimization problem. We consider the objective function
\begin{subequations}\label{eq:numeric_example}
\begin{equation}
    f(x, t_k) = \left( x_1 - e^{\cos(\omega t_k)} \right)^2 + \left( x_2 - x_1 \tanh(\sin(\omega t_k)) \right)^2,
\end{equation}
which is sampled at discrete time instances $t_k = k, \, k\in \mathbb{N}_0$. 
A slightly modified version of this function was considered in \cite{zeroing_neural}. 
One can show that this function is $m_k$-strongly convex and $L_k$-smooth with
\begin{align}
    (L_k, m_k) \in \left\{ 2 + s_k^2 \pm s_k \sqrt{s_k^2 + 4}, \, s_k = \tanh(\sin(\omega t_k)), \, k \in \mathbb{N}_0 \right\}.
\end{align}
\end{subequations}
We model this function as $f \in \Smlt$ with $\theta_{1,k} = m_k$ and $\theta_{2,k} = L_k$, and use those parameters to tune the optimization algorithm online. We consider $\omega = 0.1$, for which we have that $|m_k - m_{k-1}| \leq 0.23$ and $|L_k - L_{k-1}| \leq 0.23$.

We apply GD, 2-step GD and TM. The tracking error and the upper bounds obtained by \cref{thm:tracking:sector} and \cref{thm:tracking:variational} are visualized in Fig.~\ref{fig:tracking}. In particular, we plot a set of tracking bounds parametrized by $\rho$, where the fading colors indicate an increasing gap between bound and actual tracking error. The tracking errors confirm that TM indeed has higher tracking errors compared to the gradient descent versions, and that 2-step GD performs better than 1-step GD. For the GD algorithms, the bounds obtained by \cref{thm:tracking:sector} are closer to the actual tracking error. This can be attributed to the fact that for GD and 2-step GD, applying \cref{thm:tracking:variational} yield the same convergence rate $\rho$ but introduce additional residual terms. In contrast, for TM both bounds appear to complement each other with comparable magnitudes, where each is tighter in different timespans. Recall from Fig.~\ref{fig:iqcComparison} that the use of \cref{thm:tracking:variational} genuinely improves the convergence rate for TM, thereby offsetting the residual terms it introduces, unlike for GD and 2-step GD. This highlights that the bounds obtained by both Theorems can be both useful, and that their tightness may depend on the algorithm at hand.

\begin{figure}
  \centering
  \begin{subfigure}[t]{0.47\textwidth}
    \centering
    \input{plots/tracking/bounds_gd.tex}
    \caption{Gradient Descent.}
  \end{subfigure}
  \hfill
  \begin{subfigure}[t]{0.47\textwidth}
    \centering
    \input{plots/tracking/bounds_2gd.tex}
    \caption{2-step Gradient Descent.}
  \end{subfigure}
  \hfill
  \begin{subfigure}[t]{0.47\textwidth}
    \centering
    \input{plots/tracking/bounds_tmm.tex}
    \caption{Triple Momentum.}
  \end{subfigure}
  \caption{Tracking error of different algorithms applied to the cost \cref{eq:numeric_example} and the upper bounds obtained by \cref{thm:tracking:sector} and \cref{thm:tracking:variational}.}
  \label{fig:tracking}
  \vspace{-3mm}
\end{figure}

\section{Conclusion}
\label{sec:conclusion}

In this work, we leverage tools from robust control, LPV theory, and time-varying optimization to develop a framework to model and analyze general first-order algorithms for time-varying convex optimization. By recasting the problem as the feedback interconnection of an LPV system and a first-order oracle we are able to provide an alternative proof strategy and a computational tool to quantify its performance. For this, we introduce a variational IQC that is capable of characterizing parameter-varying nonlinearities. Using this we derive certificates for bounding the tracking error as a function of the temporal variation of the optimization problem. The use of the newly proposed variational IQCs yield better rates than the pointwise IQCs, as expected from the literature on static problems. Our bounds are novel for the proposed algorithm class, capturing several quantities associated with the temporal variability of the problem. Moreover, we show that smaller decay rates do not necessarily imply better overall tracking performance due to residual terms related to the temporal variability. This analysis framework represents a new viewpoint on time-varying optimization algorithms and can contribute to obtaining a more systematic understanding of their performance. Important open questions left for future research include robustness analysis of time-varying algorithms and the development of synthesis procedures of LPV algorithms.

\section*{Acknowledgments}
F. Jakob acknowledges the support of the International Max Planck Research School for Intelligent Systems (IMPRS-IS).
We would like to thank Prof. Peter Seiler for his valuable feedback and Dr. Lukas Schwenkel for many helpful discussions and suggestions. We further thank the Associate Editor and the anonymous reviewers, whose efforts led to
various improvements of this manuscript.

\appendix
\section{Proofs}
\label{sec:proofs}
\begin{proof}[Proof of \cref{prop:fixed_point_pstep}]
	We leverage the lifted structure of $A,\!B,\!C,\!D$ from \cref{eq:lifted_system}, fix $\theta\in\Theta$ and suppress the argument $\theta$ for ease of exposition.

	\emph{Fixed-points of \cref{eq:lifted_system}.} Note that the equilibrium relation $\xi^\star = A\xi^\star + B u^\star$ can be written equivalently as $\bar\xi_{p+1} = \bar\xi_1$, with $\bar\xi_1 := \xi^\star$ and $\bar\xi_{i+1} := \bar A\bar\xi_i + \bar B\nabla f(\bar C\bar\xi_i,\theta)$. Analogously, $y^\star = C\xi^\star + D u^\star$ unfolds to $y^{i,\star} = \bar C\bar\xi_i$. Thus the algebraic equilibrium of \cref{eq:lifted_system} is characterized by all $p$-periodic orbits of the one-step recursion. By \cref{assum:contraction}, for each fixed $\theta$, the one-step recursion is globally asymptotically stable and every trajectory converges to its unique fixed point as $i\to\infty$. Hence, no orbit of period $p>1$ can exist, and \cref{eq:lifted_system} admits exactly one fixed-point, namely $\bar\xi_i = \xi^\star$, $y^{i,\star} = \bar C \xi^\star$, $\forall i\in\mathbb{I}_p$, with $\xi^\star = \bar A\xi^\star + \bar B\nabla f(\bar C\xi^\star, \theta)$, i.e. the fixed-point of the one-step recursion \cref{eq:onestep_algo}. Note that this holds irrespective of \cref{assum:fixed_point_assum}.

	\emph{Sufficiency.} Suppose \cref{assum:fixed_point_assum} holds. From $(I-\bar A)U=0$ we obtain $\bar A^{j}U = U$ for all $j\geq 0$, and hence $AU = \bar A^{p}U = U$. Using this, we obtain $CU = [(\bar C U)^\top\ \cdots\ (\bar C\bar A^{p-1}U)^\top]^\top = 1_p\otimes(\bar C U) = 1_p\otimes I_d$. Now consider the candidate $\xi^\star = U x^\star$ and $u^\star = 0$. By $C_x U = I_d$ we have $x^\star = C_x\xi^\star$, and \cref{eq:fixpoint} becomes
	\begin{align*}
		U x^\star &= A U x^\star + B u^\star = U x^\star \\
		y^\star &= C U x^\star + D u^\star = 1_p \otimes x^\star,
	\end{align*}
	where $y^\star = 1_p\otimes x^\star$ is consistent with $u^\star = 0$ since $\nabla f(x^\star, \theta) = 0$. By the previous step, this is the unique fixed-point of \cref{eq:lifted_system}.

    \emph{Necessity.} Conversely, let \cref{eq:lifted_system} have a fixed point with $y^\star = 1_p\otimes x^\star$ and $C_x\xi^\star = x^\star$. By the first step it coincides with the fixed-point of the one-step recursion \cref{eq:onestep_algo}, whose readout is $\bar C\xi^\star = y^{1,\star} = x^\star$. Since $\theta\in\Theta$ is fixed but arbitrary and we consider all $f\in\Smlt$, the necessity part of \cref{assum:fixed_point_p1} applies and yields \cref{eq:assum:fixed_point:michalowski} as well as $\xi^\star = U x^\star$.
\end{proof}

\begin{proof}[Proof of \cref{prop:off_by_1_iqc}]

For readability, we will use the notation $g_k \triangleq \nabla f_k(x_k)$. Note that $\nabla \hat{f}_k(x_k) = \nabla f_k(x_k) = g_k$ and that
\begin{equation*}
\nabla \hat{f}_k(x_{k-1}) = \nabla \hat{f}_{k-1}(x_{k-1}) - \nabla \hat{f}_{k-1}(x_{k-1}) +  \nabla \hat{f}_k(x_{k-1}) = g_{k-1} -\Delta g_{k-1}(x_{k-1}).
\end{equation*}
Define the function
\vspace{-1mm}
\begin{equation*}
    f_k^m(x) := \hat{f}_k(x) - \frac{m_k}{2} \| x - x_k^\star \|^2
\end{equation*}
\vspace{-1mm}
and use it to define
\begin{equation*}
    V_k(x) := (L_k - m_k) f_k^m(x) - \frac{1}{2} \| \nabla f_k(x) - m_k(x - x_k^\star)\|^2.
\end{equation*}
Obviously $V_k(x_k^\star) = 0$. As shown in \cite{scherer_ZF}, $V_k(x) \geq 0$ and the inequality
\begin{equation*}\label{eq:A:proof:subgradient}
    V_k(x) - V_k(y) \leq (\nabla f_k(x) - m_k(x-x_k^\star))^\top \biggl( L_k(x-x_k^\star) - \nabla f_k(x) - (L_k(y-x_k^\star) - \nabla f_k(y)) \biggr)
\end{equation*}
holds. We plug $(x, y) \leftarrow (x_k, x_k^\star)$ to obtain
\begin{equation*}
    V_k(x_k) \leq (g_k - m_k \tilde x_k)^\top ( L_k \tilde{x}_k - g_k ).
\end{equation*} 
Moreover, we can also plug
$(x, y) \leftarrow (x_k, x_{k-1})$ to obtain
\begin{multline*}
    V_k(x_k) - V_k(x_{k-1}) \\ \leq (g_k - m_k \tilde{x}_k)^\top \biggl( L_k \tilde{x}_k - g_k - \left( L_k(\tilde{x}_{k-1} + \Delta x_{k-1}^\star) - (g_{k-1} - \Delta g_{k-1}(x_{k-1})) \right) \biggr),
\end{multline*}
where we used that $x_{k-1} - x_k^\star = \tilde{x}_{k-1} + \Delta x_{k-1}^\star$.
Now use the previous two relations to conclude
\begin{multline}
    \label{eq:appendix:subgr}
    V_k(x_k) - \rho^2 V_k(x_{k-1}) = (1-\rho^2) V_k(x_k) + \rho^2 (V_k(x_k) - V_k(x_{k-1})) \\ \leq (g_k - m_k \tilde{x}_k)^\top \biggl( L_k \tilde{x}_k - g_k - \rho^2 \left( L_k(\tilde{x}_{k-1} + \Delta x_{k-1}^\star) - (g_{k-1} - \Delta g_{k-1}(x_{k-1})) \right) \biggr).
\end{multline}
Next, for some $N\geq 1$ and $\rho \in (0,1)$, we multiply \cref{eq:appendix:subgr} with $\rho^{-2k}$ and sum from $k=0$ to $k=N$. For this, note that
\vspace{-2mm}
\begin{multline}\label{eq:appendix:lhs}
\sum_{k=0}^{N} \rho^{-2k} \left( V_k(x_k) - \rho^2 V_k(x_{k-1}) \right) \\[-2ex]
= \underbrace{\rho^{-2N}V_N(x_N)}_{\geq 0} +
 \sum_{k=1}^{N} \rho^{-2(k-1)} \left( V_{k-1}(x_{k-1}) - V_{k}(x_{k-1}) \right).
\end{multline}
\vspace{-2mm}
Plugging in the definition of $V_k$ we can expand the r.h.s. summand as 
\begin{subequations}
\begin{align}
\nonumber
& V_{k-1}(x_{k-1}) - V_{k}(x_{k-1}) \\
\label{eq:appendix:final:delta_f}
& \quad = (L_{k-1}-m_{k-1})\hat{f}_{k-1}(x_{k-1}) - (L_{k}-m_{k})\hat{f}_{k}(x_{k-1}) \\
\label{eq:appendix:final:quad:2}
& \qquad - a_{k-1}^2 \| \tilde{x}_{k-1} \|^2 + a_k^2 \| \tilde{x}_{k-1} + \Delta x_{k-1}^\star \|^2 \\
\label{eq:appendix:final:quad:3}
& \qquad - \frac{1}{2} \| g_{k-1} - m_{k-1}\tilde{x}_{k-1} \|^2 + \frac{1}{2} \| g_{k-1} - \Delta g_{k-1}(x_{k-1}) - m_{k}(\tilde{x}_{k-1}+\Delta x_{k-1}^\star) \|^2,
\end{align}
\end{subequations}
where we used the abbreviation $a_k \triangleq \sqrt{(L_k - m_k)\frac{m_k}{2}}$.
Therefore, summing $\rho^{-2k}\cref{eq:appendix:subgr}$ from $k=0$ to $k=N$, after some rearrangement, leads to
\begin{equation}
\label{eq:appendix:final}
\sum_{k=0}^N \rho^{-2k} \bigl( \left( \text{r.h.s. of } \cref{eq:appendix:subgr} \right) - \rho^2\cref{eq:appendix:final:quad:2} - \rho^2\cref{eq:appendix:final:quad:3} \bigr) \geq  \sum_{k=1}^N \rho^{-2k} \rho^2 \cref{eq:appendix:final:delta_f}.
\end{equation}

Now note that:

\noindent
(i) r.h.s. of \cref{eq:appendix:subgr}: can be realized as $(z_k^{(1)})^\top   M_1 z_k^{(1)}$ with
\begin{equation*}
    z_k^{(1)} \!=\! 
    \begin{bmatrix}
        L_k \tilde{x}_k \!-\! g_k \!-\! \rho^2  L_k(\tilde{x}_{k-1} \!+\! \Delta x_{k-1}^\star) \!+\! \rho^2 (g_{k-1} \!-\! \Delta g_{k-1}(x_{k-1})) \\
        g_k \!-\! m_k \tilde{x}_k
    \end{bmatrix}
    \;\; M_1 \!=\! \frac{1}{2}
    \begin{bmatrix}
        0 & I_d \\ I_d & 0
    \end{bmatrix}.
\end{equation*}

(ii) $-\rho^2 \cref{eq:appendix:final:quad:2}$: can be realized as $(z_k^{(2)})^\top M_2 z_k^{(2)}$ with
\begin{equation*}
    z_k^{(2)} = 
    \begin{bmatrix}
        \rho \, a_{k-1} \tilde{x}_{k-1} \\ 
        \rho \, a_{k} \left( \tilde{x}_{k-1} + \Delta x_{k-1}^\star \right)
    \end{bmatrix}
    \quad M_2 = 
    \begin{bmatrix}
        I_d & 0 \\0  & -I_d
    \end{bmatrix}.
\end{equation*}

(iii) $-\rho^2\cref{eq:appendix:final:quad:3}$: can be realized as $(z_k^{(3)})^\top M_3 z_k^{(3)}$ with
\begin{equation*}
    z_k^{(3)} = 
    \begin{bmatrix}
        \rho\left( g_{k-1} - m_{k-1}\tilde{x}_{k-1} \right) \\
        \rho \left(g_{k-1} - \Delta g_{k-1}(x_{k-1}) - m_{k}(\tilde{x}_{k-1}+\Delta x_{k-1}^\star) \right)
    \end{bmatrix}
    \quad M_3 = \frac{1}{2}
    \begin{bmatrix}
        I_d & 0 \\0  & -I_d
    \end{bmatrix}.
\end{equation*}

Combining (i)-(iii), i.e. if we define
\begin{equation*}
z_k := \left[\begin{smallmatrix}
z_k^{(1)} \\ z_k^{(2)} \\ z_k^{(3)}
\end{smallmatrix} \right], \quad \bar{M} := \mathrm{blkdiag}(M_1, M_2, M_3),
\end{equation*}
we can realize the l.h.s. of \cref{eq:appendix:final} as $\sum_{k=0}^N \rho^{-2k} z_k^\top \bar{M} z_k$. It is straightforward to check that the filter in \cref{eq:prop:off_by_1_iqc:filter} realizes such $z_k$. Therefore, \cref{eq:prop:off_by_1_iqc:ineq} corresponds to \cref{eq:appendix:final}, which concludes the proof.

\end{proof}

\bibliographystyle{siamplain}
\bibliography{references.bib}

\end{document}